\documentclass[a4paper,11pt]{amsart}
\usepackage[utf8]{inputenc}
\usepackage[T1]{fontenc}
\usepackage[english]{babel}
\usepackage{textcomp}
\usepackage{lmodern}
\usepackage{amsthm}
\usepackage{amsmath}
 \numberwithin{equation}{section}
\usepackage{amssymb}
\usepackage{mathrsfs}
\usepackage{enumerate}
\usepackage{enumitem}
\usepackage{dsfont}
\usepackage{color}
\usepackage{graphicx}
\usepackage{accents}
\usepackage[hidelinks]{hyperref}
\usepackage{float}
\usepackage{tikz, pgfplots}
\usepackage[titletoc]{appendix}
\usepackage[scale = .75]{geometry}
\newtheorem{theo}{Theorem}[section]

\newtheorem{prop}[theo]{Proposition}
\newtheorem{cor}[theo]{Corollary}
\newtheorem{lem}[theo]{Lemma}
\newtheorem{rema}[theo]{Remark}

\theoremstyle{definition}
\newtheorem{defi}[theo]{Definition}

\newcommand{\be}{\begin{equation}}
  \newcommand{\ee}{\end{equation}}
\newcommand{\bpr}{\begin{prop}}
  \newcommand{\epr}{\end{prop}}
\newcommand{\bt}{\begin{theo}}
  \newcommand{\et}{\end{theo}}
\newcommand{\bl}{\begin{lem}}
  \newcommand{\el}{\end{lem}}
\newcommand{\bc}{\begin{cor}}
  \newcommand{\ec}{\end{cor}}
\newcommand{\br}{\begin{rema}}
  \newcommand{\er}{\end{rema}}
\newcommand{\bd}{\begin{defi}}
  \newcommand{\ed}{\end{defi}}

\makeatletter
\newcommand\ds{ \displaystyle }

\newcommand\NN{\mathbb{N}}
\newcommand\RR{\mathbb{R}}

\newcommand\bU{{\mathbf U}}
\newcommand\bx{{\mathbf x}}
\newcommand\by{{\mathbf y}}
\newcommand\bn{{\mathbf n}}

\newcommand{\T}{\mathcal{T}}
\newcommand{\E}[1][]{\mathcal{E}_{#1}}

\newcommand{\Ei}[1][]{\mathcal{E}_{\text{int}#1}}
\newcommand{\Ee}[1][]{\mathcal{E}_{\text{ext}#1}}

\newcommand{\R}{\mathbb{R}}
\newcommand{\F}[1][]{F_{#1,\sigma}}

\newcommand{\Pnt}{\mathcal{P}}

\newcommand{\dx}{\,\mathrm{d}\mathbf{x}}

\newcommand{\dd}{\,\mathrm{d}}

\newcommand{\finf}{{f^\infty}}
\newcommand{\bF}{\mathbf{f}}
\newcommand{\Bpos}[1][]{B_{#1,\sigma}^{\text{pos}}}
\newcommand{\Bneg}[1][]{B_{#1,\sigma}^{\text{neg}}}

\newcommand{\meas}{{\rm m}}

\makeatother

\title{Large-time behavior of a family of finite volume schemes for boundary-driven convection-diffusion equations}

\author{Claire CHAINAIS-HILLAIRET}
\thanks{The work of Claire Chainais-Hillairet is supported by the LabEx CEMPI (ANR-11-LABX-0007-01), the MOONRISE project (ANR-14-CE23-0007) and the MoHyCon project (ANR-17-CE40-0027).}
\address{\flushleft
  Claire CHAINAIS-HILLAIRET\\[.2em]
  Univ. Lille, CNRS, UMR 8524, Inria -- Laboratoire Paul Painlevé,\\[.2em]
  F-59000 Lille, France\\[.2em]
  \texttt{E-mail: \href{mailto:claire.chainais@univ-lille.fr}{claire.chainais@univ-lille.fr}}}

\author{Maxime HERDA}
\thanks{The work of Maxime Herda is supported by a public grant overseen by the French National Research Agency (ANR) as part of the “Investissements d’Avenir” program (reference: ANR-10-LABX-0098, LabEx SMP )}
\address{\flushright
  Maxime HERDA\\[.2em]
 Inria, Univ. Lille, CNRS, UMR 8524 -- Laboratoire Paul Painlevé,\\[.2em]
  F-59000 Lille, France\\[.2em]
  \texttt{E-mail: \href{mailto:maxime.herda@inria.fr}{maxime.herda@inria.fr}}\\}

\begin{document}
    \begin{abstract}
      We are interested in the large-time behavior of solutions to finite volume discretizations of convection-diffusion equations or systems endowed with non-homogeneous Dirichlet and Neumann type boundary conditions. Our results concern various linear and nonlinear models such as Fokker-Planck equations, porous media equations or drift-diffusion systems for semiconductors. For all of these models, some relative entropy principle is satisfied and implies exponential decay to the stationary state. In this paper we show that in the framework of finite volume schemes on orthogonal meshes, a large class of two-point monotone fluxes preserve this exponential decay of the discrete solution to the discrete steady state of the scheme. This includes for instance upwind and centered convections or Scharfetter-Gummel discretizations. We illustrate our theoretical results on several numerical test cases.
    \end{abstract}
    
    \maketitle
      
      \smallskip
    

      \smallskip
    

    \tableofcontents
    
    \section{Introduction}
    
   The present paper is concerned with the study of the large-time behavior of finite volume approximations for dissipative initial-boundary problems. More precisely, we are interested in three types of convection-diffusion models that are  Fokker-Planck equations, porous media equations and drift-diffusion-Poisson systems of equations. These models are set in a bounded domain and endowed with mixed Dirichlet and no-flux boundary conditions.

   In the context of the evolution of a large number of particle, the trend to equilibrium is governed by the second law of thermodynamics. Mathematically, the dissipation of the physical entropy allows to characterize the large-time behavior through the celebrated H-theorem \cite{villani_2002_review}. Based on this principle from kinetic theory, mathematicians have intensively developed the entropy method for the study of the large time behavior of different systems of partial differential equations since the beginning of the 90's. We first refer to the survey paper by Arnold {\em et. al.} \cite{Arnold2004} for the presentation of the entropy method and for references of its application to Boltzmann equation, linear Fokker-Planck equation and  nonlinear diffusion equations for instance. Let us also mention the founding papers \cite{desvillettes_villani_2005, arnold_markowich_unterreiter_2001,carrillo_toscani_2000,carrillo_etal_2001}. Similar techniques were already used by Gajewski and Gröger in \cite{gajewski_groger_1986, gajewski_groger_1989} for the linear and nonlinear drift-diffusion systems of equations arising in semiconductor devices modeling. The entropy method has been further successfully applied to the study of the long-time behavior of reaction-diffusion equations \cite{desvillettes_fellner_2006, glitzky_gartner_2009} or cross-diffusion systems of equations \cite{burger_schlake_wolfram_2012, jungel_2015} for instance. We also refer to the recent book by Jüngel \cite{jungel_2016} and the references therein.

 In most of the literature dealing with entropy methods, models are studied without boundary conditions by considering periodic or infinite domains. In bounded domains, boundary conditions are usually chosen in a compatible way with respect to some specific equilibrium of the system. This is the case for the drift-diffusion system where the Dirichlet boundary conditions are often assumed to be compatible with the so-called thermal equilibrium. To our knowledge the first entropy method dealing with general Dirichlet-Neumann boundary conditions was proposed by  Bodineau, Lebowitz, Mouhot and Villani in \cite{bodineau_2014_lyapunov}. The key tool of our paper is the adaptation of their method to the discrete level.

 The knowledge of the large time behavior, the existence of some entropies which are dissipated along time  are  structural features, as positivity of densities or conservation of mass, that should be preserved at the discrete level by numerical schemes. The question of the large time behavior of numerical schemes has been investigated for instance for coagulation-fragmentation models \cite{filbet_2008}, for nonlinear diffusion equations \cite{chainais_jungel_schuchnigg_2016}, for reaction-diffusion systems \cite{glitzky_2008, glitzky_2011}, for drift-diffusion systems \cite{chainais_filbet_2007, bessemoulin_chainais_2017}.
 In \cite{filbet_herda_2017}, Filbet and Herda proposed a finite volume scheme preserving entropy principles for nonlinear boundary-driven Fokker-Planck equations.
 They adapted the entropy method from \cite{bodineau_2014_lyapunov} to the discrete level by defining a finite volume approximation of the steady-state and using it to design a scheme for the evolution equation.
 This pre-calculation of the steady-state guarantees that it is preserved by the evolution scheme.

 In this paper, we show that the numerical analysis performed in \cite{filbet_herda_2017} and the adaptation of the entropy method of \cite{bodineau_2014_lyapunov} also works for some usual and well-known schemes. These schemes, which include the upwind, centered and Scharfetter-Gummel scheme, do not need any pre-calculation of the steady-state, contrary to the schemes of \cite{filbet_herda_2017}. In the case of Fokker-Plank or porous media equations, we prove exponential decay towards the steady-state defined by each scheme. This constitutes the main theoretical contributions of this paper. Concerning numerical simulations, we provide several test cases illustrating these results. Moreover, we extend the scope of our theoretical results to drift-diffusion systems by some numerical tests. For this system, some special schemes (like Scharfetter-Gummel schemes), are conceived to preserve the particular form of the thermal equilibrium and ensure good long-time properties of the discrete solution. According to our simulations, all of the monotone two-point fluxes scheme tested preserve the exponential decay to their discrete equilibrium. In a nutshell, our results emphasize that discrete entropy principles and exponential decay properties may be obtained independently of the accurate approximation of the steady state.

 \subsection*{Outline}
 The outline of the paper is the following. In Section \ref{sec:fpe}, we consider linear Fokker-Planck equations with Dirichlet and no-flux boundary conditions.  We study the long-time behavior of a family of $B$-schemes (\cite{chainais_2011_finite}), including upwind, centered and Scharfetter-Gummel schemes . The main result of Section \ref{sec:fpe} is Theorem \ref{theo:fpe}: it provides a family of discrete entropy principles and shows the exponential decay towards the discrete steady-state in discrete $L^1$-norm. Section \ref{sec:pme} is devoted to porous media equations and its main result is Theorem \ref{theo:pme}. The exponential decay is established for the so-called \cite{bodineau_2014_lyapunov} relative entrophy and then for the discrete $L^{m+1}$-norm (where $m$ is the exponent in the porous medium equation). Finally, in Section \ref{sec:dds}, we consider the drift-diffusion system of equations arising in semiconductor device modeling, see for instance \cite{gajewski_groger_1986, gajewski_groger_1989}. We show by some numerical test cases that our results for Fokker-Planck equations seems to extend to these nonlinear systems. We end with some concluding remarks in Section~\ref{sec:conclusion}.

    \section{Fokker-Planck equations}\label{sec:fpe} 
    Let  $\Omega$ be a polyhedral open bounded connected subset of $\R^d$ with boundary $\Gamma = \partial\Omega$. More precisely the boundary is divided in two parts $\Gamma = \Gamma_D\cup\Gamma_N$ where $\Gamma_D$ will be endowed with non-homogeneous Dirichlet boundary conditions while $\Gamma_N$ will be endowed with no-flux boundary conditions. Let us assume that
    \[
    m(\Gamma_D)\ \neq\ 0\,,
    \]
    where $m(\cdot)$ denotes the Lebesgue measure in dimension $d-1$ (we use the same notation for the $d$-dimensional measure).
    
    The first class of models under consideration is Fokker Planck equations. It describes the evolution of a scalar density $f(t,\bx)$ by the initial-boundary value problem reading
    \begin{equation}
      \left\lbrace
      \begin{array}{rcll}
	\ds \partial_t f \,+\, {\nabla}\cdot\left({\bU}(\bx)\,f\,-\,a(\bx)\,{\nabla} f\right) &=& 0&\text{ in }\ \RR_+\times\Omega\,,\\[.75em]
	\ds f(t,\bx)&=& f^D(\bx) &\text{ on }\ \RR_+\times\Gamma_D\,,\\[.75em]
	\ds \left({\bU}(\bx)\,f\,-\,a(\bx)\,{\nabla} f\right)\cdot\bn(\bx)&=& 0 &\text{ on }\ \RR_+\times\Gamma_N\,,\\[.75em]
	\ds f(0,\bx) &=& f^\text{in}(\bx) &\text{ in }\ \Omega\,,
      \end{array}
      \right.
      \label{e:linFP}      
    \end{equation}
    The density is advected by a steady field ${\bU}:\Omega\rightarrow\R^d$ and diffused with a steady positive diffusion coefficient $a:\Omega\rightarrow\R_+$. The dynamic is driven by a non-homogeneous Dirichlet boundary condition $f^D:\Gamma_D\rightarrow\R_+$. Let us assume that $\bU\in L^\infty(\Omega)^d$, $a\in L^\infty(\Omega)$ and $f^D\in L^\infty(\Gamma_D)$ and that there are positive constants $\alpha$, $m^D$, $M^D$ and $V$ such that for almost every $\bx$ one has
    \begin{equation}
      |\bU(\bx)|\ \leq V\,,\qquad 0\,<\,\alpha\,\leq\,a(\bx)\quad\text{and}\quad 0\,<\,m^D\,\leq\,f^D(\bx)\,\leq\,M^D\,.
      \label{e:diffnondegen}
    \end{equation}

    We are interested in the long time behavior of solutions and we introduce $\finf\equiv\finf(\bx)$ a steady state of \eqref{e:linFP}, namely a function that does not depend on time and which solves the first three equations of \eqref{e:linFP}.
    Following \cite{bodineau_2014_lyapunov}, we know that there is a relative entropy structure in \eqref{e:linFP}. For any twice differentiable non-negative function $\phi$ satisfying
    \[
    \phi''>0\,,\quad\phi(1)\,=\,0\,,\quad\phi'(1)\,=\,0\,,
    \]
    there is an associated relative entropy defined by
    \be
    \mathcal{H}_\phi(t) = \int_\Omega\phi\left(\frac{f}{f^\infty}\right) f^\infty \dx\,.
    \label{e:phientropy}
    \ee
    This quantity, called relative $\phi$-entropy satisfies the following entropy/entropy-dissipation principle
    \be
    \frac{\dd \mathcal{H}_\phi}{\dd t}\,+\,\mathcal{D}_\phi\, = \,0\,,
    \label{e:entropydissipFP}
    \ee
    where the the relative $\phi$-entropy dissipation is given by
    \[
    \mathcal{D}_\phi(t) = \int_\Omega  a(\bx)\,\left|{\nabla}\left(\frac{f}{f^\infty}\right)\right|^2\,\phi''\left(\frac{f}{f^\infty}\right)\,f^\infty\, \dx\,.
    \]
    We refer to \cite[Theorem 1.4]{bodineau_2014_lyapunov} or \cite[Proposition 1.3]{filbet_herda_2017} for a proof. Typical examples of relative $\phi$-entropies are the \emph{physical relative entropy} and \emph{p-entropies} (or \emph{Tsallis relative entropies}) respectively generated, by
    \begin{equation}
      \phi_1(x)\,=\, x\ln(x) -(x-1)\,,\qquad
      \phi_p(x) \,= \,\ds\frac{x^p - px}{p-1} +1\,,\ \text{ for }\  p\in(1,2]\,.
      \label{e:examples}
    \end{equation}
    The $\phi$-entropies and their dissipations are non-negative quantities that vanish if and only if $f$ and $f^\infty$ coincide almost everywhere. The decay  of $\mathcal{H}_\phi$ characterizes the convergence to the steady state in the large. More precisely, thanks to the Poincaré inequality applied  to the function $(f-\finf)/\finf$ which vanishes on $\Gamma_D$, one can find a positive constant $\kappa$ such that the $2$-entropy and its dissipation satisfy the inequality
    \be
    \kappa\,\mathcal{H}_{\phi_2}\,\leq\mathcal{D}_{\phi_2}\,.
    \label{e:poincare}
    \ee
    As a consequence one can deduce exponential decay to equilibrium for this equation in $2$-entropy and $L^1$ norm, namely
    \[
    \|f(t)-\finf\|_{L^1(\Omega)}^2\,\leq\,\|\finf\|_{L^1(\Omega)}\mathcal{H}_{\phi_2}(t)\,\leq\,\|\finf\|_{L^1(\Omega)}\,\mathcal{H}_{\phi_2}(0)\,e^{-\kappa\,t}\,,
    \]
    where the first inequality follows from Cauchy-Schwarz inequality and the second from the combination of \eqref{e:entropydissipFP} and \eqref{e:poincare}.
    \br
    The rate $\kappa$ depends only on $\alpha$, $m^\infty$, $M^\infty$ and $\Omega$ where $m^\infty$ and $M^\infty$ are the lower and upper bounds of the unique steady state of \eqref{e:linFP} which depend additionally on $V$, $m^D$ and $M^D$. Observe that the elliptic steady equation may be non-coercive as we did not assume any sign condition on the divergence of $\bU$. Nevertheless, existence, uniqueness and the bound $M^\infty$ are shown in \cite{droniou_2002_non-coercive}. The uniform positive lower bound $m^\infty>0$ can be obtained using the De Giorgi's method \cite{degiorgi_1957, vasseur_2016_degiorgi}.
    \label{r:boundstat}
    \er
    
    In the following, we are going to show that the decay of $\phi$-entropies also holds at the discrete level for finite volume discretizations of \eqref{e:linFP} with monotone two-point flux. Let us now introduce some notation concerning these schemes. The main result of this section can be found in Theorem~\ref{theo:fpe}.

    \subsection{Numerical schemes}
    \subsubsection{Mesh and time discretization}
    Let us start with some notation associated with the discretization of $\RR_+\times\Omega$. An admissible mesh of $\Omega$ is defined by the triplet $(\T, \E, \Pnt)$. The set $\T$ is a finite family of nonempty connected open disjoint subsets $K\subset\Omega$ called control volumes or cells. The closure of the union of all control volumes is equal to the closure of $\Omega$. The set $\E$ is a finite family of nonempty subsets of $\bar{\Omega}$ called edges. Each edge is a subset of an affine hyperplane. Moreover, for any control volume $K\in\T$ there exists a subset $\E[K]$ of $\E$ such that the closure of the union of all the edges in $\E[K]$ equals to $\partial K = \bar{K}\setminus K$. We also define several subsets of $\E$. The family of interior edges is given by $\Ei = \{\sigma\in\E,\ \sigma\nsubseteq\Gamma\}$ and the family of exterior edges by $\Ee = \E\setminus\Ei$. We assume that for any edge $\sigma$, the number of control volumes sharing the edge $\sigma$ is exactly $2$ for interior edges and $1$ for exterior edges. Since every interior edge is shared by two control volumes, say $K$ and $L$, we may use the notation $\sigma = K|L$ whenever $\sigma\in\Ei$. We assume that the mesh is strongly connected, meaning that for any $K, L\in\T$ one can find a path $K_1, K_2, \dots, K_n\in\T$ such that $K_1 = K$, $K_n = L$ and $K_i$, $K_{i+1}$ share an edge for all $1\leq i\leq n$.
    
    The set $\Ee$ is the disjoint union of  $\Ee^D$ and  $\Ee^N$, composed of the exterior edges in $\Gamma_D$ and $\Gamma_N$ respectively. We make the assumption that the Dirichlet boundary condition do occur on a non-negligible subset of the boundary, namely
    \begin{equation}
      \exists\sigma\in\Ee^D\,,\quad m(\sigma)\,\neq\,0\,.
      \label{h:nonemptyEe}\tag{H1}
    \end{equation}    
    The set $\Pnt = \left\{\bx_K\right\}_{K\in\T}$ is a finite family of points satisfying that for any control volume $K\in\T$, $\bx_K\in K$. We introduce the transmissibility of the edge $\sigma$, given by 
    \[
    \tau_\sigma = \frac{m(\sigma)}{d_\sigma},
    \]
    where
    \[
    d_{\sigma} = \left\{
    \begin{aligned}
      &d(\bx_K, \bx_L),&&\text{if }\sigma\in\Ei,\ \sigma = K|L\,,\\
      &d(\bx_K, \sigma),&&\text{if }\sigma\in\Ee[,K]\,,
    \end{aligned}
    \right.
    \]
    with $d(\cdot,\cdot)$ the Euclidean distance in $\R^{d}$.
    The size of the mesh is defined by 
    \[
    \Delta x = \max_{K\in\T}\sup_{\bx,\by\in K}d(\bx,\by)\,.
    \]
    Finally we require some constraints on the mesh. For consistency of two-point gradients, we require an orthogonality condition, namely
    \begin{equation}
      \begin{aligned}
	\forall\, \bx,\by\in\sigma= K|L,\quad (\bx-\by)\cdot(\bx_K-\bx_L) = 0.
      \end{aligned}
      \label{h:orthog}
      \tag{H2}
    \end{equation}
    Additionally, in order to get discrete functional inequalities with discretization-independent constants we need the following regularity constraint on the mesh. We assume that there is a positive constant $\xi$ that does not depend on $\Delta x$ such that
    \begin{equation}
      \forall K\in\T,\ \forall \sigma\in\E[K],\quad d_{K,\sigma}\geq \xi\, d_\sigma
      \label{h:regmesh}
      \tag{H3}
    \end{equation}
    Finally, we denote the time step by $\Delta t$ and set $t^n = n\,\Delta t$.
    \subsubsection{Discrete data}
    We consider a discrete advection field  $(U_{K,\sigma})_{K\in\T,\, \sigma\in\E[K]}$, satisfying \[U_{K,\sigma} = -U_{L,\sigma}\ \text{ if }\ \sigma = K|L\,,\] as well as a discrete diffusion coefficient $(a_{\sigma})_{\sigma\in\E}$ and a discrete non-homogeneous Dirichlet boundary condition $(f_\sigma^D)_{\sigma\in\Ee^D}$. We assume that the discrete Dirichlet boundary condition is bounded from above and from below far from $0$, namely for all $\sigma\in\Ee^D$ one has 
    \be
    0\,<\,m^D\,\leq\,f_\sigma^D\,\leq\, M^D\,,
    \label{h:hyp_bc}\tag{H4}
    \ee
    for positive constants $m^D$ and $M^D$ that do not depend on the discretization. The discrete diffusion coefficient is assumed to be non-degenerate, namely there is $\alpha$ independent of the discretization such that 
    \be
    \min_{\sigma\in\E}a_\sigma \geq \alpha >0\,.
    \label{h:hyp_a}\tag{H5}
    \ee
    Moreover, we assume that the discrete advection field is bounded 
    \be
    \max_{K\in\T}\max_{\sigma\in\E[K]}|U_{K,\sigma}| \leq V\,,
    \label{h:hyp_U}\tag{H6}
    \ee
    for a constant $V$ not depending on the discretization.  
    \br
    In order to define our numerical scheme and establish our theoretical results  about long-time behavior, we do not need to specify further the discretization of the data. If one is interested in the convergence of the scheme, then if $\bU$, $a$, and $f^D$ are smooth functions, one may define their discrete counterparts as 
    \[
    U_{K,\sigma}\ =\ \frac{1}{m(\sigma)}\int_\sigma \bU(\bx)\cdot\bn_{K,\sigma}\,,\quad a_{\sigma}\ =\ \frac{1}{m(\sigma)}\int_\sigma a(\bx)\,,\quad f_{\sigma}^D\ =\ \frac{1}{m(\sigma)}\int_\sigma f^D(\bx)\,,
    \]
    where $\bn_{K,\sigma}$ is the outward normal vector of the edge $\sigma$ of the cell $K$. If the diffusion coefficient is discontinuous it is still possible to define the $a_\sigma$ in a way that will ensure consistency of the numerical fluxes defined hereafter. This requires discontinuities to be located only on edges of the mesh and, following \cite{EGHbook} and references therein, to set
    \[
    a_{\sigma}\ =\ \frac{d_\sigma\,a_K\,a_L}{d_{L,\sigma}\,a_K \,+\, d_{K,\sigma}\,a_L}\ \text{if}\ \sigma = K|L\quad \text{and}\quad a_{\sigma}\ =\ a_K\ \text{if}\ \sigma\in\Ee[,K]\,,
    \]
    with $a_K = \frac{1}{m(K)}\int_K a(\bx)$ for all cells $K\in\T$.
    \label{r:advdiff}
    \er
    \subsubsection{Definition of the numerical schemes}
    
    All the numerical schemes considered in this section can be written in the following backward Euler finite volume form
    \begin{equation}
      m(K)\frac{f_K^{n+1} - f_K^{n}}{\Delta t} \,+\, \sum_{\sigma\in\E[K]}\F[K]^{n+1}\ =\ 0\,.
      \label{e:scheme}
    \end{equation}
    At each time step $n$ one aims to find $\bF^{n+1} = (f_K^{n+1})_{K\in\T}$ from the previous $\bF^{n}$ with an initial condition $\bF^0 = (f_K^0)_{K\in\T}$.
    Now we define the flux
    \be
    \F[K]^n\ =\ \left\{
    \begin{array}{ll}
      \tau_\sigma\,a_\sigma\,\left[B(-\,U_{K,\sigma}\,d_\sigma/a_\sigma)\, f_K^n \,-\, B(U_{K,\sigma}\,d_\sigma/a_\sigma)\,\ f_{L}^n\right]&\text{ if }\sigma = K|L\\[1em]
      \tau_\sigma\,a_\sigma\,\left[B(-\,U_{K,\sigma}\,d_\sigma/a_\sigma)\, f_K^n \,-\, B(U_{K,\sigma}\,d_\sigma/a_\sigma)\,\ f_{\sigma}^D\right]&\text{ if }\sigma\in\Ee^D\\[1em]
      0&\text{ if }\sigma\in\Ee^N
    \end{array}
    \right.
    \label{e:flux}
    \ee
    The real function $B$ satisfies the following properties
    \be
    \left\{
    \begin{array}{l}
      B(x) >0\,\text{ for all }x\in\R\,\text{ and }\,B(0) \,=\, 1,\\[.75em]
      B\text{ is Lipschitz continuous}\,,\\[.75em]
      B(-x) \,-\, B(x) \,=\, x\,\text{ for all }\,x\in\R\,. 
    \end{array}
    \right.
    \label{h:hyp_B}\tag{H7}
    \ee
    These generic $B$-fluxes were introduced by Chainais-Hillairet and Droniou in \cite{chainais_2011_finite}. Observe that hypotheses \eqref{h:hyp_a}-\eqref{h:hyp_B} imply that there is $\beta$ that does not depend on the discretization such that 
    \begin{equation}
      \min_{K\in\T} \min_{\sigma\in\E[K]} B\left(\frac{|U_{K,\sigma}|d_\sigma}{a_\sigma}\right)\,\geq\, \beta\,>\,0\,.
      \label{e:condB_DeGiorgi}
    \end{equation}
    This last condition is important to prove uniform bounds on the discrete solution. 
    With particular choices of $B$ functions, one can recover several usual two-point numerical fluxes for convection-diffusion equations:
    \begin{equation}
      \begin{array}{rl}
	\text{\emph{Upwind :}}&\ds
	B(x)\ =\ 1 + \max(-x, 0),\\[1em]
	\text{\emph{Centered :}}&\ds
	B(x)\ =\ 1 - \frac{x}{2},\\[1em]
	\text{\emph{Scharfetter-Gummel :}}&\ds
	B(x)\ =\ \frac{x}{e^x-1}\text{ for all }x\neq 0\ \text{ and } B(0)=1.
      \end{array}
      \label{e:Bfunctions}
    \end{equation}
    \br
    Observe that the condition \eqref{e:condB_DeGiorgi} holds with $\beta = 1$ and without any assumption on the data in the case of the upwind scheme. 
    However, for the centered flux, this condition is not satisfied even using \eqref{h:hyp_a} and \eqref{h:hyp_U}. Actually the positivity of $B(x)$ required in \eqref{h:hyp_B} holds only if $x\,<\,2$. This issue can be resolved under a Péclet type condition on the size of the mesh. It reads
    \[
    \Delta x\ \leq\ (1-\beta)\,\frac{\alpha}{V}\,,
    \]
    for some arbitrary $\beta>0$ not depending on the discretization. This condition provides \eqref{e:condB_DeGiorgi} thanks to \eqref{h:hyp_a}, \eqref{h:hyp_U} and the fact that $d_\sigma\leq 2\Delta x$.
    \er 
    
    In order to lighten upcoming computations, let us introduce the following compact notation for the flux
    \be
    \F[K]^n\ =\ \tau_\sigma\,a_\sigma\,\left(\Bneg[K]\ f_K^n - \Bpos[K]\ f_{K,\sigma}^n\right)\,.
    \label{e:compactflux}
    \ee
    The coefficients are given by
    \[
    \left\{
    \begin{array}{lll}
      \Bneg[K]\ =\ B(-\,U_{K,\sigma}\,d_\sigma/a_\sigma)\,,&\Bpos[K]\ =\ B(+\,U_{K,\sigma}\,d_\sigma/a_\sigma)&\text{if }\sigma\in\E\setminus\Ee^N\,,\\[.5em]
      \Bneg[K]\ =\ \Bpos[K]\ =\ 0&&\text{if }\sigma\in\Ee^N\,.
    \end{array}\right.
    \]
    They are always non-negative thanks to \eqref{h:hyp_B}. The neighbor unknown is  
    \be
    f_{K,\sigma}^n\ =\ \left\{
    \begin{array}{rcl}
      f_L^n&\text{ if }&\sigma = K|L\,, \\[.5em]
      f_\sigma^D&\text{ if }&\sigma\in\Ee^D\,,\\[.5em]
      f_K^n&\text{ if }&\sigma\in\Ee^N\,.
    \end{array}
    \right.
    \label{e:discretebc}
    \ee  
    and we additionally define the difference operator
    \be
        D_{K,\sigma} f^n\ =\ f_{K,\sigma}^n\,-\,f_K^n\,,\ \forall K\in\T,\forall \sigma \in\E[K].
    \label{e:diffop}
    \ee
    One has the relation 
    \[
    \Bneg[K]\ =\ \Bpos[L]\,,
    \]
    so that it is clear that the flux defined in \eqref{e:flux_rewritten} satisfies the conservativity property $F_{K,\sigma}\,+\,F_{L,\sigma}\,=\,0$ for interior edges $\sigma = K|L$.
    
    \subsection{Discrete steady state}
    
    Let us now turn to the corresponding stationary version of the problem. We say that $\bF^{\infty} = (f_K^\infty)_{K\in\T}$ is a discrete steady state of the scheme \eqref{e:scheme}-\eqref{e:flux} if
    \be
    \sum_{\sigma\in\E[K]}\F[K]^\infty\ =\ 0\,,\quad\text{ for all }K\in\T\,.
    \label{e:steady_divfree}
    \ee
    The steady flux is defined by
    \be
    \F[K]^\infty\ =\ \left\{
    \begin{array}{ll}
      \tau_\sigma\,a_\sigma\,\left[B(-\,U_{K,\sigma}\,d_\sigma/a_\sigma)\, f_K^\infty \,-\, B(U_{K,\sigma}\,d_\sigma/a_\sigma)\,\ f_{L}^\infty\right]&\text{ if }\sigma = K|L\\[1em]
      \tau_\sigma\,a_\sigma\,\left[B(-\,U_{K,\sigma}\,d_\sigma/a_\sigma)\, f_K^\infty \,-\, B(U_{K,\sigma}\,d_\sigma/a_\sigma)\,\ f_{\sigma}^D\right]&\text{ if }\sigma\in\Ee^D\\[1em]
      0&\text{ if }\sigma\in\Ee^N
    \end{array}
    \right.
    \label{e:steadyflux}
    \ee
    which may be rewritten
    \[
    \F[K]^\infty\ =\  \tau_\sigma\,a_\sigma\,\left(\Bneg[K]\ f_K^\infty - \Bpos[K]\ f_{K,\sigma}^\infty\right)\,,
    \]
    with the same notation as in \eqref{e:compactflux}.
    
    Let us rewrite this in matrix form as $\mathbb{M}\,\mathbf{f}^\infty = \mathbf{b}^D$. The matrix $\mathbb{M}$ has the following coefficients
    \[
    M_{K,K}\ =\ \sum_{\sigma\in\E[K]}\tau_\sigma\,a_\sigma\,\Bneg[K]\quad\text{and}\quad
    M_{K,L}\ =\ \left\{
    \begin{array}{ll}
      -\tau_\sigma\,a_\sigma\,\Bpos[K]&\text{ if there is }\ \sigma\ =\ K|L\,,\\[.5em]
      0&\text{ otherwise}\,.
    \end{array}
    \right.
    \]
    and the vector $\mathbf{b}^D = (b^D_K)_{K\in\T}$ contains the values at the boundary, that is \[b^D_K\, =\, \sum_{\sigma\in\E[K]\cap\Ee^D}\tau_\sigma\,a_\sigma\,\Bpos[K]\,f_\sigma^D\,.\] 
    \bpr
    Under hypotheses \eqref{h:nonemptyEe}-\eqref{h:hyp_U} there exists a unique discrete steady state $(f_K^\infty)_{K\in\T}$ satisfying the steady scheme  \eqref{e:steady_divfree}. Moreover there are positive constants $m^\infty, M^\infty$ depending only on $\Omega$, $m^D$, $M^D$, $\alpha$, $V$ and $\xi$ such that
    \[
    0\,<\,m^\infty\,\leq\,f_K^\infty\,\leq\, M^\infty\,.
    \]
    for all $K\in\T$.
    \label{p:exist_steady}
    \epr
    \begin{proof}
      
      In order to prove existence, uniqueness and non-negativity of the discrete solution it suffices to prove that $\mathbb{M}$ is a non-singular M-matrix (see \cite[Chapter 6, Definition 1.2]{Berman_Plemmons_1994}). We first show that $\mathbb{M}$ is an M-matrix and then that it is non-singular. 
      Observe that the coefficients of $\mathbb{M}$ satisfy
      \[
      |M_{K,K}| \,-\,\sum_{\sigma = K|L}|M_{L,K}|\ =\ \sum_{\sigma\in\E[K]\cap\Ee^D}\tau_\sigma\,a_\sigma\,\Bneg[K]\ \geq\  0\,.
      \]
      Hence $\mathbb{M}$ is diagonally dominant by columns. Therefore for any $\varepsilon>0$, if we denote by $\mathbb{I}$ the identity matrix, $\mathbb{M} + \varepsilon \mathbb{I}$ is strictly diagonally dominant. By \cite[Chapter 6, Theorem 4.6, (C9)]{Berman_Plemmons_1994}, since $\mathbb{M} + \varepsilon \mathbb{I}$ is non-singular for any $\varepsilon>0$,  $\mathbb{M}$ is an M-matrix. Now observe that
      \[
      |M_{K,K}| \,-\,\sum_{\sigma = K|L}|M_{L,K}|\ >\ 0\quad\text{ if }\E[K]\cap\Ee^D\neq\emptyset\,.
      \]
      Since the mesh is strongly connected and $\Ee^D\neq\emptyset$, $\mathbb{M}$ is a diagonally dominant matrix with a nonzero elements chain (see \cite{shivakumar_1974_sufficient}). From the main theorem of \cite{shivakumar_1974_sufficient}, $\mathbb{M}$ is invertible.
      
      The uniform bounds on the solution can be obtained by adapting the De Giorgi method \cite{degiorgi_1957}, \cite{vasseur_2016} to the discrete setting. For the lower bound, we can adapt the result of Chainais-Hillairet, Merlet and Vasseur in \cite{chainais_merlet_vasseur_2017} with minor modifications. The important assumptions are the $L^\infty$ bound \eqref{h:hyp_U} and the bound from below \eqref{e:condB_DeGiorgi}. Concerning the upper bound $M^\infty$, a similar strategy can be adopted. The main difference is to replace truncated solution by its composition with the function $x\mapsto\ln(1+x)$ in order to initialize the iterative argument. We refer to \cite{droniou_2002_non-coercive} for a proof in this spirit at the continuous level.   
    \end{proof}
    
    \br[Incompressible advections] 
    Assume the whole boundary is endowed with Dirichlet boundary conditions, that is $\Ee^N = \emptyset$. Then if the advection field is divergence free, namely if
    \[
    \sum_{\sigma\in\E[K]}m(\sigma)\,U_{K,\sigma}\ =\ 0\quad\text{ for all }K\in\T\,,
    \]
    then a strong minimum and maximum principle hold, namely $m^\infty = m^D$ and $M^\infty = M^D$. Indeed, under the previous assumptions one has for all $K\in\T$ that 
    \[
    \sum_{L\in\T}M_{K,L}(f_L^\infty - m^D)\ =\ \sum_{\sigma\in\E[K]\cap\Ee^D}\tau_\sigma\,a_\sigma\,\Bpos[K](f_\sigma-m^D)\,\geq\,0\,,
    \]
    so $\mathbb{M}\,(\mathbf{f^\infty} - m^D) \geq 0$ and similarly $\mathbb{M}\,(M^D -\mathbf{f^\infty})\geq 0$ (with component-wise substraction and inequalities). By the monotony of $\mathbb{M}$, one has $(\mathbf{f^\infty} - m^D) \geq 0$ and $(M^D -\mathbf{f}^\infty)\geq 0$. When the assumptions are not satisfied there is potentially overshoot ($M^\infty \geq M^D$) and undershoot ($m^\infty \leq m^D$).
    \er

    \subsection{Reformulation using the steady flux}
    The key observation in our analysis is that we can reformulate the flux at time $t^n$ using the discrete steady state and flux. Indeed, 
    the flux \eqref{e:compactflux} can be rewritten in the following ways
    \[     
    \begin{array}{rcccl}
      \F[K]^n &=& \F[K]^\infty\, h_K^n&-& \tau_\sigma\,a_\sigma\,\Bpos[K]\,f_{K,\sigma}^\infty\,D_{K,\sigma}h^n\\[.5em]
      &=& \F[K]^\infty\, h_{K,\sigma}^n &-& \tau_\sigma\,a_\sigma\,\Bneg[K]\,f_{K}^\infty\,\,D_{K,\sigma}h^n\,.
    \end{array}
    \]
    where
    \[
    h_K^n\ =\ \frac{f_K^n}{f_K^\infty}\ \text{ for all }K\in\T\quad \text{and}\quad  h^D_\sigma\ =\ 1\text{ for }\sigma\in\Ee^D\,.
    \]
    From there we can rewrite the flux in the following upwind form $(h_K^n)_{K\in\T}$
    \be
    \F[K]^n\ =\ \left(\F[K]^\infty\right)^+ h_K^n \,-\, \left(\F[K]^\infty\right)^- h_{K,\sigma}^n  \,-\, \tau_\sigma\,a_\sigma\,f_{B,\sigma}^{\infty}\,D_{K,\sigma}h^n\,.
    \label{e:flux_rewritten}
    \ee
    where $u^+ = \max(u,0)$ and $u^- = \max(-u,0)$ and the discrete steady state on edges is defined by 
    \be
    f_{B,\sigma}^{\infty}\ =\ \min(\Bneg[K]\,f_K^\infty,\ \Bpos[K]\,f_{K,\sigma}^\infty).
    \label{e:discrete_edge_var}
    \ee
    \br Observe that $f_{B,\sigma}^{\infty}$ is well-defined as it does not depend on $K$ since $\Bpos[K] = \Bneg[L]$ when $\sigma = K|L$. Moreover, since $B(0) = 1$, this quantity is consistent with the continuous steady state thanks to \eqref{h:hyp_U} and \eqref{h:hyp_B}. Finally, on exterior edges, beware that $f_{B,\sigma}^{\infty}$ need not be equal to $f_\sigma^{\infty}$.
    \er
    Once the flux has been rewritten as \eqref{e:flux_rewritten}, it enters the framework introduced by Filbet and Herda in \cite{filbet_herda_2017} and entropy dissipation properties can be established.
    \subsection{Entropy dissipation and long-time behavior}
    
    Now we define the discrete counterpart of relative $\phi$-entropies and their dissipations. As in the continuous setting we consider a twice continuously differentiable function $\phi$ satisfying $\phi'(1) = \phi(1) = 0$ and we define the discrete relative $\phi$-entropy at time $t^n$ by 
    \be
    H_\phi^n\ =\ \sum_{K\in\T} m(K)\,\phi\left(h_K^n\right)\,f^\infty_K\,.
    \label{e:ent_phi}
    \ee
    and the discrete dissipations of $\phi$-relative entropy at time $t^n$ as
    \[
    D_\phi^n\ =\ \sum_{\sigma\in\E}\tau_\sigma\,a_\sigma\,(D_{K,\sigma}h^n)\,(D_{K,\sigma}\phi'(h^n))\, f_{B,\sigma}^{\infty}\,\geq\,0\,,
    \]
    where the non-negativity stems from the monotony of $\phi'$. The main result of this section is the following.
    
    \bt  Let us choose any function $B$, advection field, diffusion coefficient and boundary conditions such that \eqref{h:nonemptyEe}-\eqref{h:hyp_U} hold. Then for any non-negative initial condition $(f^0_K)_{K\in\T}$, there is a unique solution to the scheme 
    \eqref{e:scheme}-\eqref{e:flux}. Moreover, it satisfies
    \begin{equation}
      \frac{H_{\phi}^{n+1} - H_{\phi}^{n}}{\Delta t}\  +\ D_\phi^{n+1} \leq 0\,, \quad\text{ for all }n\in\NN\,.
      \label{e:decr_phi_imp}
    \end{equation}
    Besides, the discrete solution $(f_K^n)_{K\in\T,n\in\NN}$ satisfies the uniform bounds
    \[
    m^\infty\min\left(1,\min_{K\in\T}f^0_K/f^\infty_K\right)\ \leq\, f_K^n\ \leq\, M^\infty\,\max\left(1,\max_{K\in\T}f^0_K/f^\infty_K\right)\,.
    \]
    uniformly for  $K\in\T$ and $n\in\NN$,
    Finally the discrete solution decays exponentially fast in time to the discrete steady state of the scheme $(f_K^\infty)_{K\in\T}$ in the following sense. For any $k >0$, if $\Delta t \leq k$  there is a constant $\kappa$ depending only on $\Omega$, $\xi$, $m^\infty$, $M^\infty$, $\alpha$, $\beta$ and $k$ such that for all $n\in\NN$ one has 
    \[
    H_{\phi_2}^{n}\ \leq \ H_{\phi_2}^{0}\,\,e^{-\kappa\, t^n}\,,
    \]
    and 
    \[
    \left(\sum_{K\in\T}m(K)|f_K^n-f_K^\infty|\right)^{2}\leq \ H_{\phi_2}^{0}\,\left(\sum_{K\in\T}m(K)|f_K^\infty|\right)\,e^{-\kappa\, t^n}.
    \]
    \label{theo:fpe}
    
    \et
    \br 
    The decay rate is at least 
    \[
    \kappa = \frac{1}{k}\ln\left(1+k\,\xi\,\beta\,\frac{\alpha\,m^\infty}{C_P\,M^\infty}\right)
    \]
    where $C_P$ depends only on the domain. The expected rate at the continuous level is $\frac{\alpha\,m^\infty}{C_P\,M^\infty}$.
    \er
    \begin{proof} We proceed in four steps.
	
	\emph{Step 1: Existence, uniqueness and non-negativity.} With the same notations as in Proposition~\ref{p:exist_steady}, one can reformulate the scheme as $(\mathbb{I} + \Delta t \mathbb{M})\mathbf{f}^{n+1}\ =\ \mathbf{f}^n + \Delta t\mathbf{b}^D$, where $\mathbb{I}$ is the identity matrix and $\mathbf{f}^n$ the vector of unknown at time $t^n$. Since $\mathbb{M}$ is an M-matrix, so is  $\mathbb{I} + \Delta t \mathbb{M}$. Hence the scheme has a unique solution and it is non-negative.
	
	
	\emph{Step 2: Entropy / entropy dissipation inequality.}  By convexity of $\phi$ one has
	\[
	\begin{array}{rcl}
	  H^{n+1}_\phi - H^{n}_\phi &\leq&\ds \sum_{K\in\T}m(K)(f_K^{n+1} - f_K^{n})\phi'(h^{n+1}_K)\\[1em]
	  &\leq&\ds-\Delta t\,\sum_{K\in\T}\sum_{\sigma\in\E[K]}F_{K,\sigma}^{n+1}\,\phi'(h^{n+1}_K)
	\end{array}
	\]
	where we used the scheme \eqref{e:scheme}. Now let us show that 
	\be
	\sum_{K\in\T}\sum_{\sigma\in\E[K]}\F[K]^{n} \,\phi'\left(h^{n}_K\right)\ \geq\  D_\phi^{n}\,.
	\label{e:discrete_dissip}
	\ee
	We use the reformulated flux \eqref{e:flux_rewritten}
	and treat separately the convection and diffusion part. The flux writes
	\[
	\F[K]^n\ =\ \F[K]^\text{upw}  \,-\, \tau_\sigma\,a_\sigma\,f_{B,\sigma}^{\infty}\,D_{K,\sigma}h^n\,,
	\]
	with the convective upwind part
	\[
	\F[K]^\text{upw}\ =\ \left(\F[K]^\infty\right)^+ h_K^n \,-\, \left(\F[K]^\infty\right)^- h_{K,\sigma}^n\,.
	\]
	The diffusion part provides the dissipation with an integration by parts  
	\begin{multline*}
	  -\sum_{K\in\T}\sum_{\sigma\in\E[K]}\tau_\sigma\,a_\sigma\,f_{B,\sigma}^{\infty}\,(D_{K,\sigma}h^n) \,\phi'\left(h_K^n\right)\\\ =\ \sum_{\sigma\in\E}\tau_\sigma\,a_\sigma\,f_{B,\sigma}^{\infty}\,(D_{K,\sigma}h^n)\,(D_{K,\sigma}\phi'(h^n))\ =\ D_\phi^n\,,
	\end{multline*}
	It remains to prove that the contribution of $\F[K]^\text{upw}$ in \eqref{e:discrete_dissip} is non-negative. For this, we define the $\phi$-centered flux
	\[
	\F[K]^\text{cen}\ =\ \F[K]^\infty\, M^\phi(h_K^n,h_{K,\sigma}^n)\,,
	\]
	where the $\phi$-mean \cite{filbet_herda_2017} is defined by
	\[
	M^\phi(s,t)\ =\ \frac{\varphi(s) - \varphi(t)}{\phi'(s)-\phi'(t)},
	\]
	with  $\varphi(s) = s\phi'(s) - \phi(s)$. Observe that the function $M^\phi$ is symmetric, satisfies $\min(s,t)\leq M^\phi(s,t)\leq\max(s,t)$ and $M^\phi(s,s) = s$. Now let us integrate in the discrete sense the difference between the two fluxes against $\phi'(h)$, namely
	\[
	\begin{array}{rcl}
	  &&\ds\sum_{K\in\T}\sum_{\sigma\in\E[K]}\left(\F[K]^\text{upw} - \F[K]^\text{cen}\right) \,\phi'\left(h_K^n\right)\\[.75em]
	  &=&\ds -\sum_{\sigma\in\E}\left(\F[K]^\text{upw} - \F[K]^\text{cen}\right) \,(D_{K,\sigma}\phi'(h^n))\\[.75em]
	  &=&\ds \sum_{\sigma\in\E}\left(\F[K]^\infty\right)^+ \left(M^\phi(h_K^n,h_{K,\sigma}^n) - h_{K}^n\right)\,(D_{K,\sigma}\phi'(h^n))\\[.75em]
	  &&+ \ds\sum_{\sigma\in\E}\left(\F[K]^\infty\right)^- \left(h_{K,\sigma}^n - M^\phi(h_K^n,h_{K,\sigma}^n)\right)\,(D_{K,\sigma}\phi'(h^n))
	  \ \geq\  0\,.
	\end{array}
	\]
	The non-negativity is a consequence of the monotony of $\phi'$ and of the properties of $M^\phi$. Finally observe that by definition of $M^\phi$ one has
	\[
	\begin{array}{rcl}
	  \ds\sum_{K\in\T}\sum_{\sigma\in\E[K]}\F[K]^\text{cen} \,\phi'\left(h_K\right)
	  &=&\ds -\sum_{\sigma\in\E}\F[K]^\infty\, M^\phi(h_K,h_{K,\sigma})\,(D_{K,\sigma}\phi'(h^n))\\[.75em]
	  &=&\ds -\sum_{\sigma\in\E}\F[K]^\infty\, (D_{K,\sigma}\varphi'(h^n))\\[.75em]
	  &=& \ds\sum_{K\in\T}\varphi'\left(h_K\right)\,\sum_{\sigma\in\E[K]}\F[K]^\infty\\[.75em]&=& 0\,,
	\end{array}
	\]
	where we used the definition of the steady state \eqref{e:steady_divfree} in the last equality. Eventually, combining the last two computations, one gets that
	\[
	\sum_{K\in\T}\sum_{\sigma\in\E[K]}\F[K]^\text{upw} \,\phi'\left(h_K\right)\ \geq\ 0\,.
	\]
	It follows that \eqref{e:discrete_dissip} holds and consequently so does \eqref{e:decr_phi_imp}.
	
	\emph{Step 3: Uniform upper and lower bounds.} From the entropy inequality \eqref{e:decr_phi_imp}, one has in particular that $0\leq H_\phi^n\leq H_\phi^0$ for all admissible functions $\phi$. By an approximation argument these inequalities hold also for the convex functions $\phi_{+}(x) = (x-M)^+$ and $\phi_{-}(x) = (x-m)^-$ for any $M\geq1\geq m$. Take $M = \max\left(1,\max_{K\in\T}h^0_K\right)$ and $m=\min\left(1,\min_{K\in\T}h^0_K\right)$. Then $H_{\phi_+}^0 = H_{\phi_-}^0 =0$ and therefore $H_{\phi_+}^n = H_{\phi_-}^n =0$ which yields the uniform bounds on the discrete solution.
	
	\emph{Step 4: Exponential decay.} We use the discrete Poincaré inequality (see \cite{EGHbook} or \cite[Theorem 6]{bessemoulin_chainais_filbet_2015}), one has 
	\[
	\begin{array}{rcl}
	  H_{\phi_2}^n&\leq&\ds M_\infty\,\sum_{K\in\T}m(K)\left(h_K^n-1\right)^2\\[1em]
	  &\leq&\ds M_\infty\,\frac{C_P}{\xi}\sum_{\sigma\in\E}\tau_\sigma\left(D_{K,\sigma}(h^n-1)\right)^2\\[1em]
	  &\leq&\ds\frac{C_P\,M_\infty}{\beta\,\xi\,\alpha\,m_\infty}D_{\phi_2}^n\,.
	\end{array}
	\]
	Combining this with the $\phi_2$-entropy inequality yields 
	\[
	H_{\phi_2}^{n+1}\,\leq\,\left(1 + \frac{\beta\,\xi\,\alpha\,m_\infty}{C_P\,M_\infty}\,\Delta t\right)^{-1}H_{\phi_2}^{n}\,,
	\]
	for all $n\in\NN$. It provides the exponential decay of $H_{\phi_2}^{n}$. Exponential decay in $L^1$ is then a consequence of the Cauchy-Schwarz inequality.
    \end{proof}
    \br
    The estimates obtained in Theorem~\ref{theo:fpe} provide uniform discrete $L^2(\RR_+;\,H^1)$ and $L^\infty$ control on the discrete solution. This is enough to prove the convergence of approximate solutions to weak solutions of \eqref{e:linFP} if the data is discretized as in Remark~\ref{r:advdiff}. The proof follows classical arguments that the interested reader may find in \cite{EGHbook} and \cite{chainais_2011_finite}.
    \label{r:convergence}
    \er

    \subsection{Numerical results}\label{sec:numresFP}
\begin{figure}[!h]
%
%
\begin{tikzpicture}

\begin{axis}[%
width=.3\linewidth,
height=.3\linewidth,
at={(0.758in,0.499in)},
scale only axis,
unbounded coords=jump,
xmin=0,
xmax=1,
xlabel style={font=\color{white!15!black}},
xlabel={$x_1$},
ymin=0,
ymax=1,
ylabel style={font=\color{white!15!black}},
ylabel={$x_2$},
axis background/.style={fill=white}
]
\addplot [color=black, forget plot]
  table[row sep=crcr]{%
0	0.5\\
0.25	0.5\\
nan	nan\\
0	0.5\\
0	0.75\\
nan	nan\\
0	0.5\\
0.15	0.65\\
nan	nan\\
0	0.5\\
0	0.25\\
nan	nan\\
0	0.5\\
0.175	0.325\\
nan	nan\\
0.25	0.5\\
0.5	0.5\\
nan	nan\\
0.25	0.5\\
0.15	0.65\\
nan	nan\\
0.25	0.5\\
0.325	0.675\\
nan	nan\\
0.25	0.5\\
0.35	0.35\\
nan	nan\\
0.25	0.5\\
0.175	0.325\\
nan	nan\\
0.5	0.5\\
0.5	0.75\\
nan	nan\\
0.5	0.5\\
0.325	0.675\\
nan	nan\\
0.5	0.5\\
0.75	0.5\\
nan	nan\\
0.5	0.5\\
0.65	0.65\\
nan	nan\\
0.5	0.5\\
0.5	0.25\\
nan	nan\\
0.5	0.5\\
0.35	0.35\\
nan	nan\\
0.5	0.5\\
0.675	0.325\\
nan	nan\\
0.5	0.75\\
0.5	1\\
nan	nan\\
0.5	0.75\\
0.325	0.675\\
nan	nan\\
0.5	0.75\\
0.35	0.85\\
nan	nan\\
0.5	0.75\\
0.65	0.65\\
nan	nan\\
0.5	0.75\\
0.675	0.825\\
nan	nan\\
0.5	1\\
0.25	1\\
nan	nan\\
0.5	1\\
0.35	0.85\\
nan	nan\\
0.5	1\\
0.75	1\\
nan	nan\\
0.5	1\\
0.675	0.825\\
nan	nan\\
0.25	1\\
0	1\\
nan	nan\\
0.25	1\\
0.35	0.85\\
nan	nan\\
0.25	1\\
0.175	0.825\\
nan	nan\\
0	1\\
0	0.75\\
nan	nan\\
0	1\\
0.175	0.825\\
nan	nan\\
0	0.75\\
0.15	0.65\\
nan	nan\\
0	0.75\\
0.175	0.825\\
nan	nan\\
0.15	0.65\\
0.325	0.675\\
nan	nan\\
0.15	0.65\\
0.175	0.825\\
nan	nan\\
0.325	0.675\\
0.35	0.85\\
nan	nan\\
0.325	0.675\\
0.175	0.825\\
nan	nan\\
0.35	0.85\\
0.175	0.825\\
nan	nan\\
0.75	0.5\\
1	0.5\\
nan	nan\\
0.75	0.5\\
0.65	0.65\\
nan	nan\\
0.75	0.5\\
0.825	0.675\\
nan	nan\\
0.75	0.5\\
0.85	0.35\\
nan	nan\\
0.75	0.5\\
0.675	0.325\\
nan	nan\\
1	0.5\\
1	0.75\\
nan	nan\\
1	0.5\\
0.825	0.675\\
nan	nan\\
1	0.5\\
1	0.25\\
nan	nan\\
1	0.5\\
0.85	0.35\\
nan	nan\\
1	0.75\\
1	1\\
nan	nan\\
1	0.75\\
0.825	0.675\\
nan	nan\\
1	0.75\\
0.85	0.85\\
nan	nan\\
1	1\\
0.75	1\\
nan	nan\\
1	1\\
0.85	0.85\\
nan	nan\\
0.75	1\\
0.85	0.85\\
nan	nan\\
0.75	1\\
0.675	0.825\\
nan	nan\\
0.65	0.65\\
0.825	0.675\\
nan	nan\\
0.65	0.65\\
0.675	0.825\\
nan	nan\\
0.825	0.675\\
0.85	0.85\\
nan	nan\\
0.825	0.675\\
0.675	0.825\\
nan	nan\\
0.85	0.85\\
0.675	0.825\\
nan	nan\\
0	0\\
0.25	0\\
nan	nan\\
0	0\\
0	0.25\\
nan	nan\\
0	0\\
0.15	0.15\\
nan	nan\\
0.25	0\\
0.5	0\\
nan	nan\\
0.25	0\\
0.15	0.15\\
nan	nan\\
0.25	0\\
0.325	0.175\\
nan	nan\\
0.5	0\\
0.5	0.25\\
nan	nan\\
0.5	0\\
0.325	0.175\\
nan	nan\\
0.5	0\\
0.75	0\\
nan	nan\\
0.5	0\\
0.65	0.15\\
nan	nan\\
0.5	0.25\\
0.325	0.175\\
nan	nan\\
0.5	0.25\\
0.35	0.35\\
nan	nan\\
0.5	0.25\\
0.65	0.15\\
nan	nan\\
0.5	0.25\\
0.675	0.325\\
nan	nan\\
0	0.25\\
0.15	0.15\\
nan	nan\\
0	0.25\\
0.175	0.325\\
nan	nan\\
0.15	0.15\\
0.325	0.175\\
nan	nan\\
0.15	0.15\\
0.175	0.325\\
nan	nan\\
0.325	0.175\\
0.35	0.35\\
nan	nan\\
0.325	0.175\\
0.175	0.325\\
nan	nan\\
0.35	0.35\\
0.175	0.325\\
nan	nan\\
0.75	0\\
1	0\\
nan	nan\\
0.75	0\\
0.65	0.15\\
nan	nan\\
0.75	0\\
0.825	0.175\\
nan	nan\\
1	0\\
1	0.25\\
nan	nan\\
1	0\\
0.825	0.175\\
nan	nan\\
1	0.25\\
0.825	0.175\\
nan	nan\\
1	0.25\\
0.85	0.35\\
nan	nan\\
0.65	0.15\\
0.825	0.175\\
nan	nan\\
0.65	0.15\\
0.675	0.325\\
nan	nan\\
0.825	0.175\\
0.85	0.35\\
nan	nan\\
0.825	0.175\\
0.675	0.325\\
nan	nan\\
0.85	0.35\\
0.675	0.325\\
nan	nan\\
};
\end{axis}
\end{tikzpicture}%
      \caption{Coarsest mesh ($\Delta x = 1/4$).}\label{f:mesh}
    \end{figure}
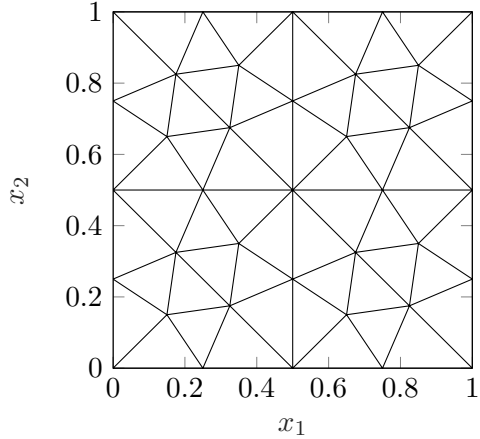

    Let us now illustrate these theoretical results on two test cases. In the following the domain is set to $\Omega = [0,1]\times[0,1]$ and the simulations are performed on a family of meshes generated from the one shown on Figure~\ref{f:mesh}. The refinement of a mesh is obtained taking $2$ by $2$ grids of its contraction by a $1/2$ factor.
    
    \subsubsection{Toy-model}
    The first test case (inspired from \cite{filbet_herda_2017}) is the following 2D Fokker-Planck equation
    \[
    \partial_t f \,+\, {\nabla}\cdot\left(\bU\,f\,-\,{\nabla} f\right)\ =\ 0\,. 
    \]
    with advection given by the vector $\bU=\binom{1}{0}$ and endowed with the boundary conditions $f(0,x_2) = 1$, $f(1,x_2) = \exp(1)$ on the left and right edges $\Gamma_D = (\{0\}\times[0,1])\cup(\{1\}\times[0,1]$) and no-flux boundary conditions on the top and bottom edges $\Gamma_N = ([0,1]\times\{0\})\cup([0,1]\times\{1\}$). An explicit solution of this equation is given by  
    \be\label{explicit}
    f(t,x_1,x_2)\ =\ \exp(x_1) + \exp\left(\frac{x_1}{2}-\left(\pi^2+\frac 14\right)t\right)\sin(\pi x_1)\,.
    \ee
    The corresponding steady state is $\finf(x_1,x_2)\ =\ \exp(x_1)$. 
    
    We first solve the steady equation on our meshes using the finite volume scheme \eqref{e:steady_divfree}-\eqref{e:steadyflux} for each function $B$ defined in \eqref{e:Bfunctions}. Observe that the advection field derives from a potential, namely $\bU = \nabla\varphi$ where $\varphi(x_1,x_2) = x_1$. Following this expression we define the discrete advection as follows $U_{K,\sigma}\,d_\sigma\ =\ D_{K,\sigma}\varphi$, where $\varphi_K\ =\ \varphi(\bx_K)$ for each control volume $K\in\T$ and for $\sigma\in\Ee$, $\varphi_\sigma\ =\ \varphi(\bx_\sigma)$ with $\bx_\sigma$ the orthogonal projection of $\bx_K$ on $\sigma$. With this discretization one expects the upwind and centered schemes to provide a discretization of the steady state of order $1$ and $2$ respectively. The Scharfetter-Gummel scheme enjoys the nice property of being exact for the steady equation since $\finf = \exp(\phi)$. By exact, we mean that $(\exp(\phi(\bx_K)))_{K\in\T}$ solves the discrete Scharfetter-Gummel scheme and we shall refer to this steady state as the \emph{``real discrete steady state''} (as opposed to \emph{``the steady state of the scheme''}).
    
    Experimentally, we compute the $L^1$ error between the \emph{steady state of each scheme} and \emph{real discrete steady state}. The results are given in Table~\ref{tb:order1global} and confirms the previous claims. Observe that the Scharfetter-Gummel scheme is indeed exact and numerical errors deteriorates while refining the mesh due to the addition of machine epsilons.

    \begin{table}[!h]
      \[
      \begin{array}{|c||c|c||c|c||c|}
	\hline	
	&\textbf{Upwind} &&\textbf{Centered}&&\textbf{Scharfetter-}\\
	& &&&&\textbf{Gummel}\\
	\Delta x&\text{Error in}\ L^1&\text{Order}&\text{Error in}\ L^1&\text{Order}&\text{Error in}\ L^1\\
	\hline
	1/4&6.04.10^{-3}&&1.23.10^{-4}&&9.15.10^{-16}\\
	1/8&3.24.10^{-3}&0.90&3.05.10^{-5}&2.01&3.96.10^{-15}\\
	1/16&1.67.10^{-3}&0.95&7.67.10^{-6}&1.99&4.50.10^{-15}\\
	1/32&8.50.10^{-4}&0.98&1.92.10^{-6}&1.99&5.05.10^{-15}\\
	1/64&4.28.10^{-4}&0.99&4.83.10^{-7}&2.00&2.19.10^{-14}\\
	1/128&2.15.10^{-4} &0.99&1.21.10^{-7}&2.00&2.36.10^{-14}\\
	1/256&1.08.10^{-4}&1.00&3.03.10^{-8}&2.00&1.18.10^{-13}\\
	\hline
      \end{array}
      \]   
      \caption{Error in $L^1$ and experimental order of convergence between the \emph{steady state of the scheme} and the \emph{real steady state}. }\label{tb:order1global}
    \end{table}

    \begin{figure}[!h]
      \includegraphics[width = .8\textwidth]{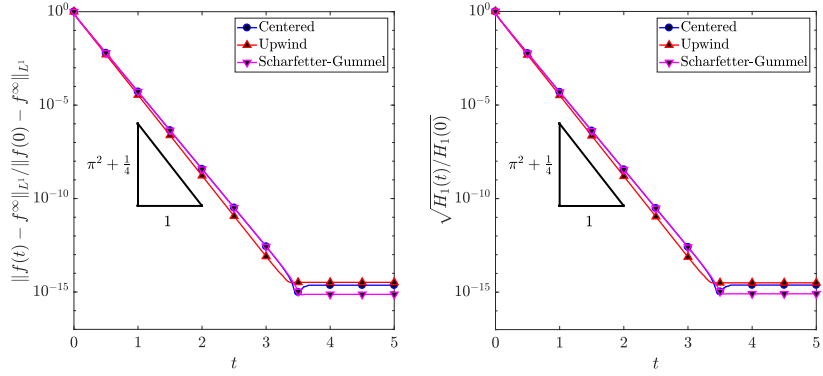}
      \caption{Decay to the \emph{steady state associated to each scheme} in $L^1$ norm (left) and Boltzmann entropy (right) for each scheme. Here the mesh is that of Figure~\ref{f:mesh} ($\Delta x = 1/4$). Time step: $\Delta t = 10^{-2}$.}\label{f:ownsteady}
    \end{figure}

    \begin{figure}[!h]
      \includegraphics[width = .8\textwidth]{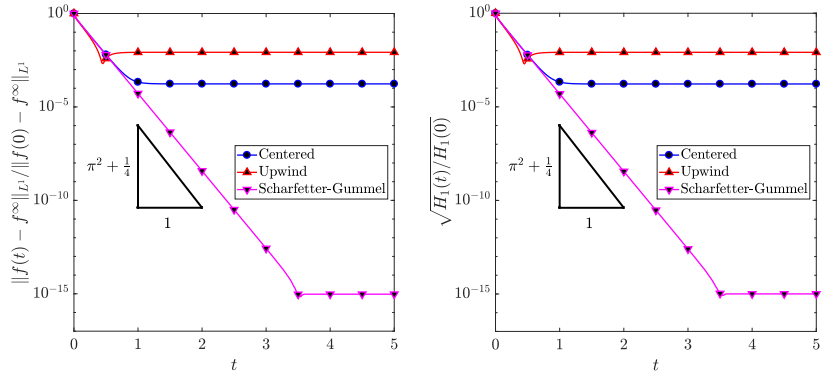}
      \caption{Decay to the \emph{real discrete steady state} in $L^1$ norm (left) and Boltzmann entropy (right) for each scheme.}\label{f:realsteady}
    \end{figure}
    
    Secondly, we turn to the resolution of the evolution equations and to the long-time behavior of solutions. We use the scheme \eqref{e:scheme}-\eqref{e:flux} for each function $B$ defined in \eqref{e:Bfunctions} and with the discretization of the convection described above. Analytically, from the explicit formula \eqref{explicit} one expects the solution to decay to the steady state with exponential rate $\pi^2+1/4$ both in $L^1$ and in square roots of $\phi$-entropies. From the result of Theorem~\ref{theo:fpe}, at the discrete level, we expect exponential decay to the \emph{steady state associated to each scheme}. This is illustrated on Figure~\ref{f:ownsteady} and is in accordance with the theoretical result. On Figure~\ref{f:realsteady} we see that experimentally this exponential decay do not hold for every scheme when the solution is compared to the real discrete state. Indeed, convergence reaches a threshold before machine precision for both the upwind and centered schemes. The comparison between Figure~\ref{f:ownsteady} and Figure~\ref{f:realsteady} emphasizes that the accurate long-time behavior of a numerical scheme is relative to the definition of the discrete steady state.

    \subsubsection{Diffusion in an heterogeneous media with gravity}
\begin{figure}[!h]
      \centering
      \includegraphics[width = \linewidth]{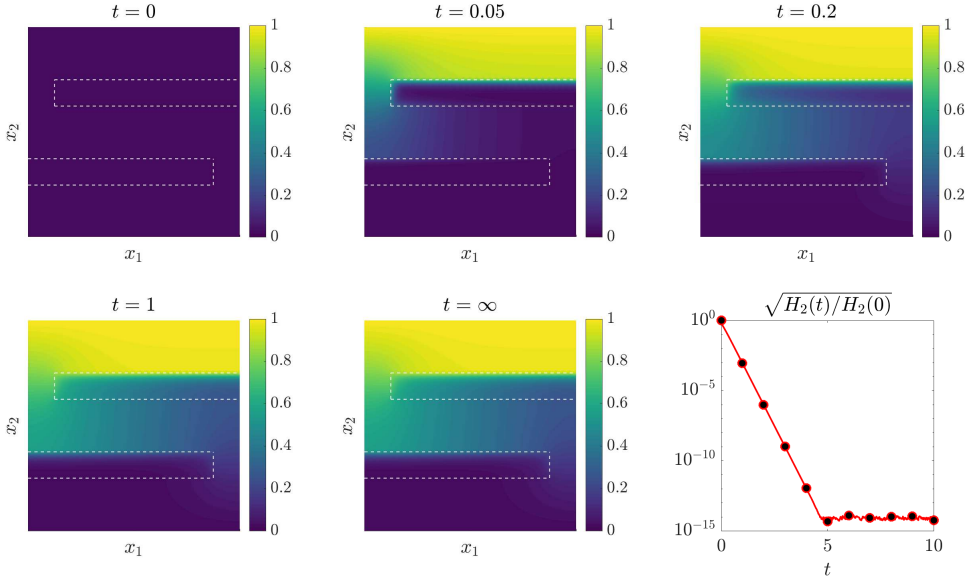}\\
      \caption{Snapshots of the solution at different time, snapshot of the steady state of the scheme and decay to steady state in $2$-entropy. Here the mesh is that of Figure~\ref{f:mesh} after four refinements ($\Delta x = 1/64$). Time step: $\Delta t = 10^{-2}$. The upwind scheme was used for this simulation.}\label{f:cas_test_cances_guichard}
    \end{figure}

    The following test case is inspired by the fourth test case of \cite{cances_guichard_2017}. We consider an heterogeneous convection-diffusion 
    \[
    \partial_t f \,+\, {\nabla}\cdot\left(\bU\,f - a(\bx)\,{\nabla} f\right)\ =\ 0\,. 
    \]
    More precisely, the domain $\Omega = [0,1]^2$ is separated in two open subdomains $\Omega_1$ and $\Omega_2$ which represents a drain and a barrier respectively. On the snapshots of Figure~\ref{f:cas_test_cances_guichard}, $\Omega_2$ is the union the two open subsets surrounded by the white dashed line and $\Omega_1$ is the interior of $\Omega\setminus\Omega_2$. The diffusion coefficient takes the following values 
    \[
    a(\bx)\ =\ \left\{
    \begin{array}{ll}
      a_1\ =\ 3&\text{if }\bx\in\Omega_1\,,\\[.5em]
      a_2\ =\ 0.01&\text{if }\bx\in\Omega_2\,.
    \end{array}
    \right.
    \]
    The fourth mesh is used so that the interface between $\Omega_1$ and $\Omega_2$ is entirely made of cell edges. The diffusion coefficient is discretized as explained in Remark~\ref{r:advdiff}. The initial data is the constant function $f_0(\bx) = 0.018$ for all $\bx\in\Omega$. The boundary conditions are of Dirichlet type on the top and bottom boundaries and of no-flux type on the left and right boundaries. The boundary data is $f^D(x_1,x_2) = 1$ on the top boundary $x_2 = 1$ and $f^D(x_1,x_2) = 0.018$ on the bottom boundary $x_2 = 0$.
    The advection field given by the vector $\bU=\binom{-1/2}{0}$.
    The results of Figure~\ref{f:cas_test_cances_guichard} illustrate the exponential decay of relative entropy of the scheme, as predicted by Theorem~\ref{theo:fpe}.

    \section{Porous medium equations}\label{sec:pme}
    
    In this section, we consider the porous medium equation with Dirichlet-Neumann boundary conditions.
    
    Let $f^D\in L^{\infty}({\Gamma^D})$, $f^{\rm in}\in L^{\infty}(\Omega)$ satisfying \eqref{e:diffnondegen} and $f^{\rm in}\geq 0$ and let $m> 1$. The porous medium equation  is a nonlinear scalar diffusion equation reading 
    \begin{equation}
      \left\lbrace
      \begin{array}{ll}
	\ds \partial_t f \,=\, \Delta f^m&\text{ for }\ \bx\in\Omega, \  t\geq0\,,\\[.75em]
	\ds f(t,\bx)\,=\, f^D(\bx) &\text{ for }\ \bx\in\Gamma_D, \  t\geq0\,,\\[.75em]
	\ds {\nabla} f\cdot\bn(\bx)\,=\, 0 &\text{ for }\ \bx\in\Gamma_N, \  t\geq0\,,\\[.75em]
	\ds f(0,\bx) \,=\, f^\text{in}(\bx) &\text{ for }\ \bx\in\Omega\,,
      \end{array}
      \right.
      \label{e:porous}      
    \end{equation}
    For existence and uniqueness results concerning the porous medium equation \eqref{e:porous}, we refer to the monograph by V\'azquez \cite{vazquez2007}. As in Section~\ref{sec:fpe}, we are interested in the long time behavior of solutions. We consider $\finf$ the steady state associated to \eqref{e:porous} and defined by
    \[
      \left\lbrace
      \begin{array}{ll}
	\Delta(\finf)^m=0&\text{ for }\ \bx\in\Omega,\\[.75em]
	\ds \finf(\bx)\,=\, f^D(\bx) &\text{ for }\ \bx\in\Gamma_D, \\[.75em]
	\ds {\nabla} \finf\cdot\bn(\bx)\,=\, 0 &\text{ for }\ \bx\in\Gamma_N,
      \end{array}
      \right.
      \]	
    
    There is a wide literature concerning the large-time asymptotics of the porous medium equation. The pionneering work by Alikakos and Rostamian \cite{alikakos_rostamian_81} deals with the case of homogeneous Neumann boundary conditions: the sharp decay rate $t^{-1/(m-1)}$ in the $L^\infty$ norm is established for general initial data and an exponential decay rate is also shown for strictly positive initial data. Similar results are proved by Grillo and Muratori in \cite{grillo_muratori_2013, grillo_muratori_2014} using weighted Sobolev inequalities or Gagliardo-Nirenberg inequalities.
    The case of homogeneous Dirichlet boundary conditions is studied in \cite{ grillo_muratori_2014}.  In \cite{chainais_jungel_schuchnigg_2016}, 
    Chainais-Hillairet, J\"{u}ngel and Schuchnigg also establish algebraic or exponential decay of zeroth-order and first-order entropies for the porous-medium and the fast-diffusion ($m<1$) equations with homogeneous Neumann boundary conditions. Their proof is based on Beckner-type functional inequalities and extends to the discrete case, so that they get the long-time asymptotics of finite volume approximations of the porous-medium and the fast-diffusion equations.
    
    In order to study the long time behavior of the porous medium equation with general Dirichlet-Neumann boundary conditions, Bodineau, Lebowitz, Mouhot and Villani introduce in \cite{bodineau_2014_lyapunov} the relative entrophy $\mathcal{N}_m$ defined by      
    \be
    \mathcal{N}_m(t)\ =\ \int_\Omega\frac{f^{m+1}-(\finf)^{m+1}}{m+1}-(\finf)^m(f-\finf)\dx.
    \label{entropyPM}
    \ee
    Observe that with the notations \eqref{e:examples}, the entropy can be reformulated as 
    \[
    \mathcal{N}_m(t) = \frac{m+1}{m}\,\int_\Omega\phi_{m+1}\left(\frac{f}{\finf}\right)\,(\finf)^{m}\,\dx\,,
    \]
    which does not coincide with the ${m+1}$-entropy \eqref{e:phientropy}. A direct computation shows that the entrophy is associated to the dissipation 
    \[
    \mathcal{D}_m(t)\ =\ \int_\Omega\left|\nabla(f^m-\finf^m)\right|^2\,\dx,
    \]
    through the relation 
    \be
    \frac{\dd \mathcal{N}_m}{\dd t}\,+\,\mathcal{D}_m\,=\,0\,.
    \label{e:entropydissipPM}
    \ee
    Thanks to the Poincaré inequality and an elementary functional inequality (which we will detail in what follows), $\mathcal{N}_m$ is controlled by $\mathcal{D}_m$, which ensures the exponential decay of $\mathcal{N}_m$ and the convergence $f\to  \finf$  as $t\to\infty$ (\cite[Theorem 1.8]{bodineau_2014_lyapunov}). 
    The aim is now to adapt this strategy to the discrete level, in order to show the long time behavior of a classical finite-volume approximation  for the porous medium equation with general Dirichlet-Neumann boundary conditions.
    
    \subsection{Numerical scheme}
    
    We consider the classical finite volume scheme with two-point flux approximation for the porous medium equation and its steady state. 
    The notations for the mesh are the same as in Section \ref{sec:fpe}. For $g=(g_K)_{K\in\T}$ and $g^D=(g_\sigma^D)_{\sigma\in\E[{\rm ext}]^D}$, we have already defined in  \eqref{e:discretebc} and \eqref{e:diffop} the quantities $g_{K,\sigma}$  and $ D_{K,\sigma} g$ for all  for all $K\in \T$, and $\sigma \in \E[K]$.     Let us remark that the definition of $g_{K,\sigma}$ in \eqref{e:discretebc} ensures that $D_{K,\sigma} g=0$ for all $\sigma\in \E[{\rm ext}]^N$.
    
    The scheme for the porous medium equation is a backward Euler scheme in time. It writes
    \begin{equation}\label{scheme:pme}
      \meas (K)\frac{f_K^{n+1}-f_K^n}{\Delta t}\,+\,\sum_{\sigma\in\E[K]} \tau_\sigma D_{K,\sigma} (f^{n+1})^m\ =\ 0\quad \forall K\in\T
    \end{equation}
    and starts from a given discretization of the initial data $(f_K^0)_{K\in\T}$. On the boundary we consider a given discretization $(f_\sigma^D)_{\sigma\in\Ee^D}$ of the non-homogeneous Dirichlet condition and we refer to Remark~\ref{r:advdiff} for an example of such discretization. We assume that the boundary condition satisfies the uniform upper and lower bounds \eqref{h:hyp_bc}. 
    Existence and uniqueness of a solution to the scheme \eqref{scheme:pme} has been established by Eymard, Gallouët, Herbin and Nait Slimane in \cite[Remark 2.3]{eymard_naitslimane_98} (see also \cite{EGHbook}). 
    
    The scheme for the steady-state writes :
    \begin{equation}\label{scheme:pme_stat}
      \sum_{\sigma\in\E[K]} \tau_\sigma D_{K,\sigma} (\finf)^m\ =\ 0\,,\quad \forall K\in\T\,.
    \end{equation}
    Setting $u_K=(f_K^\infty)^m$ for all $K\in\T$ with $u_\sigma^D=(f_\sigma^D)^m$ for all $\sigma \in \E[{\rm ext}]^D$, the scheme 
    \eqref{scheme:pme_stat} rewrites as the classical finite volume scheme for the Laplace equation, which admits a unique solution if \eqref{h:nonemptyEe} holds. 
    \br[Neumann boundary conditions]
    Uniqueness is lost if $\Gamma^D=\emptyset$. However, in this case, the steady state to \eqref{scheme:pme} has necessarily the same mass as the initial condition:
    \[
    \ds\sum_{K\in\T} \meas(K) f_K^\infty =\ds\sum_{K\in\T} \meas(K) f_K^0, 
    \]
    so that 
    \[
    f_K^\infty=\ds\frac{1}{\meas(\Omega)}\sum_{K\in\T} \meas(K) f_K^0\,,
    \]
    for all $K\in\T$.
    \label{r:fullneumann}
    \er

    \subsection{Exponential decay of the entrophy}

    Let $(f_K^n)_{K\in\T, n\geq 0}$ be the solution to the scheme \eqref{scheme:pme} and $(f_K^\infty)_{K\in\T}$ be the solution to the scheme \eqref{scheme:pme_stat}. We define the discrete relative entrophy by 
    \begin{equation}\label{Nmdis}
      \mathcal{N}_m^n=\sum_{K\in\T}\meas(K)\left(\frac{(f_K^{n})^{m+1}-(f_K^{\infty})^{m+1}}{m+1}-(f_K^{\infty})^m(f_K^{n}-f_K^\infty)\,.
      \right)\quad \forall n\geq 0.
    \end{equation}
    The discrete dissipation associated to \eqref{Nmdis} is given by
    \[
      \mathcal{D}_m^{n+1}=\ds\sum_{\sigma \in \E}\tau_{\sigma} \left( D_{K,\sigma}((f^{n+1})^m-(f^{\infty})^m)\right)^2\,.
    \]

    The following theorem is the main result of Section~\ref{sec:pme}. It states the exponential decay with respect to time of the discrete entrophy \eqref{Nmdis}. As a consequence, it also provides the exponential decay of the $L^{m+1}$-norm.
    
    \begin{theo}\label{theo:pme}
      
      Under hypotheses \eqref{h:nonemptyEe}-\eqref{h:hyp_bc}, the solutions $(f_K^n)_{K\in\T, n\geq 0}$ and $(f_K^\infty)_{K\in\T}$ respectively to the scheme \eqref{scheme:pme} and the scheme \eqref{scheme:pme_stat} are such that 
      \be\label{dNdt_dis:pme}
      \ds\frac{\mathcal{N}_m^{n+1}-\mathcal{N}_m^n}{\Delta t}\, +\, \mathcal{D}_m^{n+1}\ \leq\ 0\,,\quad  \forall n\in\NN\,.
      \ee
      As a consequence for any $k>0$ and $\Delta t\leq k$ there is a positive constant $\lambda$ depending only on $k$, $\Gamma^D$, $\Omega$, $\xi$, $m^D$ and $m$, so that 
      \be\label{expdecay_dis:pme}
      \mathcal{N}_m^n\ \leq\ e^{-\lambda t^n} \mathcal{N}_m^0\,,\quad \forall n\in\NN\,.
      \ee
      and 
      \be\label{expdecay2_dis:pme}
      \ds\sum_{K\in\T} \meas(K)|f_K^{n}-f_K^{\infty}|^{m+1}\ \leq\ (m+1)\,e^{-\lambda t^n} \mathcal{N}_m^0\,,\quad \forall n\in\NN\,.
      \ee
    \end{theo}
    \br\label{r:rate}
    The decay rate can be taken uniformly as 
    \[
    \lambda = \frac{1}{k}\ln\left(1+k\,\xi\,\frac{{(m^D)}^{m-1}}{C_P}\right)
    \]
    where $C_P$ depends only on the domain. The expected rate at the continuous level is ${{(m^D)}^{m-1}}/{C_P}$.
    \er
    \begin{proof} 
      The proof can be split into 3 steps. We first establish the discrete counterpart of the entropy / entropy dissipation relation \eqref{e:entropydissipPM}. Then, thanks to a discrete Poincar\'e inequality and an elementary inequality, we rely the discrete dissipation to the entropy, which leads to \eqref{dNdt_dis:pme} and \eqref{expdecay_dis:pme}. Finally, we obtain \eqref{expdecay2_dis:pme} thanks to another elementary inequality.
	
\emph{Step 1: Entrophy dissipation.} Let us first prove a discrete counterpart of \eqref{e:entropydissipPM}. We have
	\[
	\mathcal{N}_m^{n+1}- \mathcal{N}_m^n=\sum_{K\in\T}\meas(K)\left(\frac{(f_K^{n+1})^{m+1}-(f_K^{n})^{m+1}}{m+1}-(f_K^{\infty})^m(f_K^{n+1}-f_K^n) 
	\right).
	\]
	But the convexity of the function $x\mapsto x^{m+1}$ implies that $y^{m+1}-x^{m+1}\leq (m+1) y^m (y-x)$ for all $x,y\in \R$, so that 
	\[
	\mathcal{N}_m^{n+1}- \mathcal{N}_m^n \leq \sum_{K\in\T}\meas(K)\left((f_K^{n+1})^{m}-(f_K^{\infty})^{m}\right)(f_K^{n+1}-f_K^n).
	\]
	Using the schemes \eqref{scheme:pme} and \eqref{scheme:pme_stat} and applying a discrete integration by parts, we obtain 
	\[
	\mathcal{N}_m^{n+1}- \mathcal{N}_m^n\ \leq\ -\Delta t\,\mathcal{D}_m^{n+1}\,.
	\]
	
\emph{Step 2:  Entrophy / Entrophy dissipation inequality.} Let us now recall the discrete Poincar\'e inequality (see \cite{EGHbook} or \cite[Theorem 6]{bessemoulin_chainais_filbet_2015}). Thanks to \eqref{h:nonemptyEe} one has
	\[
	\ds\sum_{K\in\T} \meas (K) (g_K)^2\ \leq\ \frac{C_P}{\xi} \sum_{\sigma\in\E} \tau_{\sigma} (D_{K,\sigma} g)^2,
	\]
	for every set of discrete values $(g_K)_{K\in\T}$ associated to zero boundary values $g_{\sigma}=0$ for all $\sigma\in\Ee^D$. We may apply this inequality to $g= (f^{n+1})^m-(f^{\infty})^m$. It yields
	\[
	\ds\sum_{K\in\T} \meas (K) ((f_K^{n+1})^m-(f_K^{\infty})^m)^2\ \leq\ \frac{C_P}{\xi} \mathcal{D}_m^{n+1}.
	\]
	But, for $m\geq 1$ we have the following elementary inequality (whose proof is left to the reader)
	\[
	(z^m-1)^2\ \geq\ \ds\frac{1}{m+1}(z^{m+1}-(m+1)z+m)\,,\quad\forall z\geq 0.
	\]
	As a consequence, we get that 
	\[
	\left((f_K^{n+1})^m-(f_K^{\infty})^m\right)^2\ \geq\ (f_K^{\infty})^{m-1}\left(\ds\frac{(f_K^{n+1})^{m+1}-(f_K^{\infty})^{m+1}}{m+1}-(f_K^{\infty})^m(f_K^{n+1}-f_K^\infty)\right)    
	\]
	and 
	\be \label{linkND:pme}
	{\mathcal{N}}_{m}^{n+1} \leq \frac{C_P}{\xi\,{(m^D)}^{m-1}} {\mathcal{D}}_m^{n+1}\,,\quad \forall n\in\NN\,.
	\ee
	Then, we deduce from \eqref{dNdt_dis:pme} and \eqref{linkND:pme} that
	\[
	\mathcal{N}_{m}^{n+1}\ \leq\ \left(1\,+\,\Delta t\,\frac{\xi\,{(m^D)}^{m-1}}{C_P}\right)^{-1}\,\mathcal{N}_{m}^{n}\,,\quad \forall n\in\NN\,.
	\]
	It yields \eqref{expdecay_dis:pme} as a straightforward consequence.
	
\emph{Step 3: Decay in $L^{m+1}$-norm.} It remains to prove \eqref{expdecay2_dis:pme}. It is just a consequence of the following inequality (once again left to the reader)
	$$
	\vert z-1\vert^{m+1}\leq z^{m+1}-(m+1)z+m\,,\quad \forall z\geq 0\,,
	$$
	which implies that 
	$$
	\ds\sum_{K\in\T} \meas (K) \vert f_K^n-f_K^\infty\vert^{m+1}\leq (m+1) {\mathcal{N}}_m^n\,, \quad \forall  n\in\NN\,,
	$$
	and leads to \eqref{expdecay2_dis:pme} thanks to \eqref{expdecay_dis:pme}.
    \end{proof}

    \subsection{Numerical results}
    
    In the following, two test cases are presented to illustrate our theoretical results. The domain and the mesh structure are the same as in Section~\ref{sec:numresFP}.
    
    \subsubsection{Implementation}\label{implementation}
     The scheme for the steady state \eqref{scheme:pme_stat} is a linear system of equations on the new unknown $\mathbf{u} = ((f_K^\infty)^m)_{K\in\T}$. However, whenever the exponent $m > 1$, the numerical scheme for the evolution equation \eqref{scheme:pme} is nonlinear since it is implicit in time. Hence one has to solve a nonlinear system of equations at each time step. In practice, we use a Newton method and an adaptive time step strategy in order to approximate the solution. At time step $n+1$, we launch a Newton method starting at $\bF^n$ to solve the equation. If the method has not converged with precision $\varepsilon = 10^{-11}$ before $50$ steps, we divide the time step by $2$ and restart the Newton method. At the beginning of each time step we multiply the previous time step by $2$. The initial time step is $\Delta t_0 = 10^{-3}$ and we impose additionally that $10^{-8} \leq \Delta t_n \leq 10^{-2}$ to avoid over-refinement or coarsening of the time step during the procedure.

    \subsubsection{Filling of a porous media}
    
    \begin{figure}[!h]
      \centering
      \includegraphics[width = .8\linewidth]{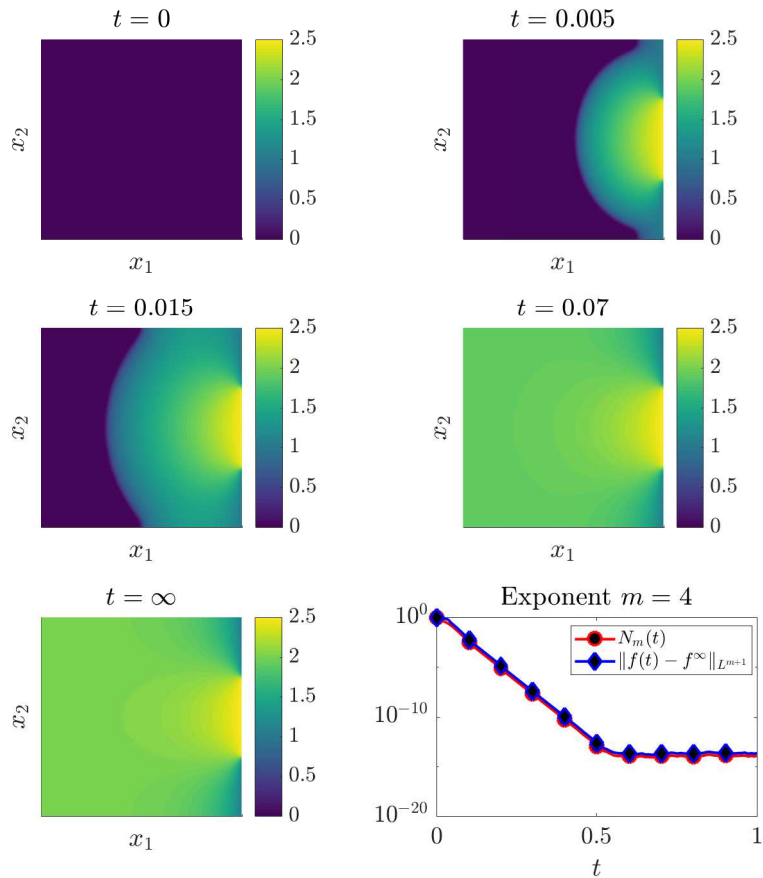}\\
      \caption{Snapshots of the solution at different times, snapshot of the steady state of the scheme and decay to steady state in entrophy and $L^{m+1}$-norm.  $\Delta x\ =\ 1/64$. }\label{f:cas_test_porous}
    \end{figure}
    
    In this first test case, we illustrate the exponential decay to equilibrium in entrophy and $L^{m+1}$ norm. The porous medium equation is taken with exponent $m = 4$. The initial data is $f^\text{in}(\bx) = 0$ for $\bx\in\Omega$. The right edge $\Gamma_D = {1}\times[0,1]$ of the square domain is endowed with Dirichlet boundary conditions $f_D(\bx) = 2.5$ if $x_1 = 1$ and $x_2\in(0.3, 0.7)$, and  $f_D(\bx) = 1$ if $x_1 = 1$ and $x_2\in[0, 0.3]\cup[0.7, 1]$. The rest of the edge $\Gamma_N = \Gamma\setminus\Gamma_D$ is endowed with Neumann boundary conditions. The results are displayed on Figure~\ref{f:cas_test_porous}. One can see qualitatively the finite propagation speed in the media as well as the exponential decay in time predicted by Theorem~\ref{theo:pme}.

    \subsubsection{Influence of parameters on the decay rate}

    \begin{figure}[!h]
      \includegraphics[width = .8\linewidth]{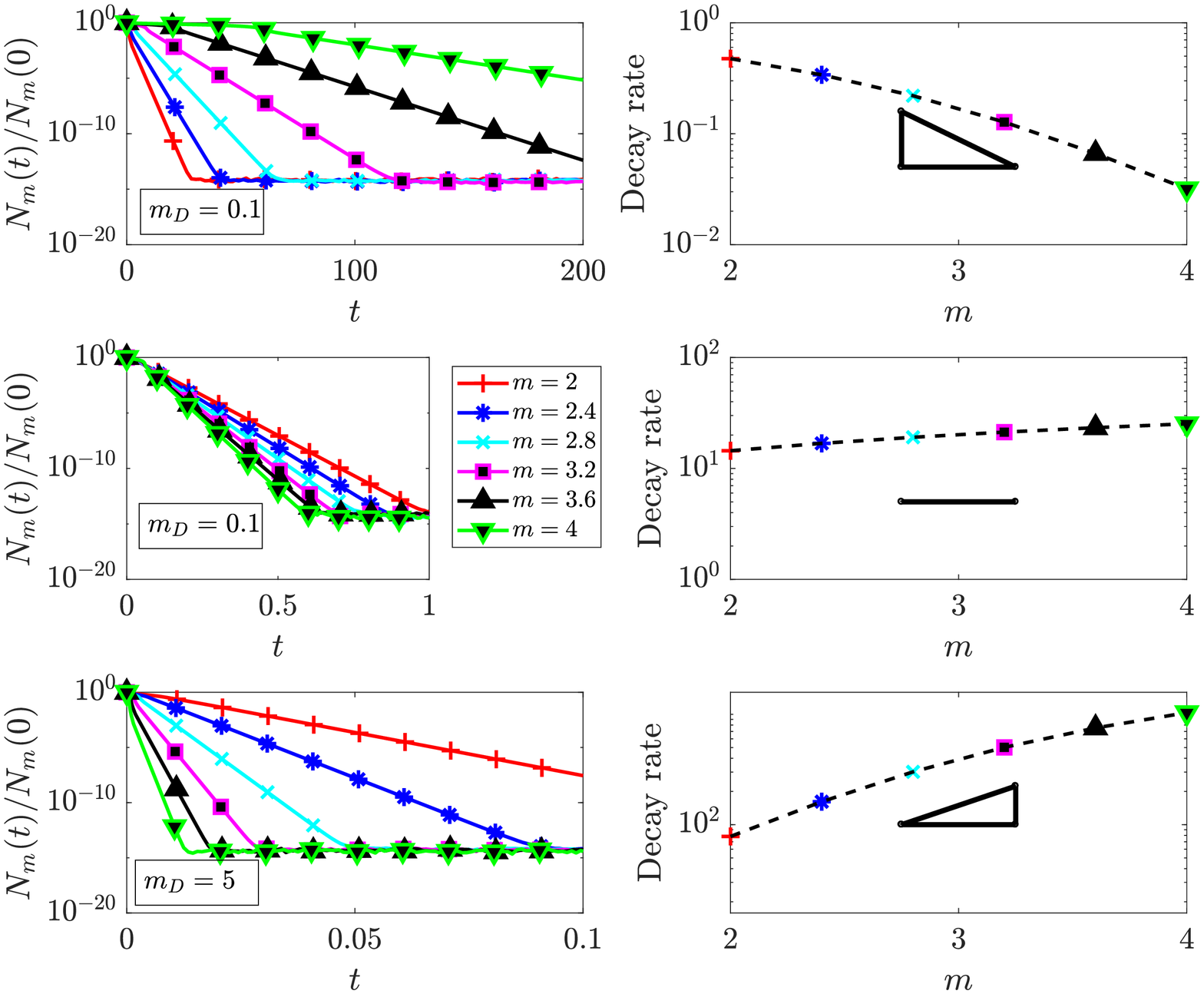}\\
      \caption{Dependence of the decay rates in $m$ and $m_D$. \emph{Left:} $L^{m+1}$-norm versus $t$ for different values of $m$ and $m_D$. \emph{Right:} experimental decay rate versus $m$. The triangles show the predicted slope from Remark~\ref{r:rate}.   $\Delta x\ =\ 1/64$. }\label{f:cas_test_porous_parameters}
    \end{figure}
    
    In this second test case we want to illustrate numerically the influence of the parameters of the model on the decay rate and compare it with the rate announced in Remark~\ref{r:rate}. More precisely, we vary the exponent $m \in\{2,\,2.4,\,2.8,\,3.2,\,3.6,\,4\}$ and the minimal Dirichlet boundary value $m_D\in\{0.1,\,1,\,5\}$. The initial data  is  $f^\text{in}(\bx) = 0$ for $\bx\in\Omega$. The whole boundary is endowed with Dirichlet boundary conditions $f_D(\bx) = m_D$.  The results are displayed on Figure~\ref{f:cas_test_porous_parameters}. They show that the experimental rate behaves nearly as $O(m_D^{m-1})$ in terms of monotony and linear dependence between the logarithm of the rate and $m$.
    
    \section{Drift-diffusion systems}\label{sec:dds}
    The third model under consideration is the Van Roosebroeck's hydrodynamic equations for semiconductors. It is a nonlinear drift-diffusion system, consisting of two continuity equations controlling time evolutions of the electron density  $N(t,\bx)$ and the hole density $P(t,\bx)$ coupled with the Poisson equation for the electrostatic potential $V(t,\bx)$. It reads 
    \be\label{eq:NPV}
    \left\{
    \begin{array}{rcl}
      \partial_tN \, - \, \nabla\cdot\left(\mu_N\,(\nabla N - N\,\nabla V\right)) & = & 0\,,\\[.5em]
      \partial_tP \, - \, \nabla\cdot\left(\mu_P\,(\nabla P + P\,\nabla V\right)) & = & 0\,,\\[.5em]
      -\lambda^2\,\Delta V\ =\ P\,-\,N\,+\,C\,,&&
    \end{array}
    \right.
    \ee
    where $C(\bx)$ is a background doping profile characterizing the device. The dimensionless physical parameters $\mu_{N}$, $\mu_{P}$ and $\lambda$ are respectively the mobilities of electrons and holes, and the Debye length. In the following, we take $\mu_N = \mu_P = 1$. The system \eqref{eq:NPV} is supplemented with the initial conditions $N_0$, $P_0$ 
    and two types of boundary conditions. The first part of the boundary $\Gamma_D$ is endowed with the non-homogeneous Dirichlet boundary conditions  denoted by $N^D$, $P^D$, $V^D$    corresponding physically to ohmic contacts. The rest of the boundary $\Gamma_N$ is insulated with homogeneous Neumann boundary conditions.    
    
    Existence and uniqueness of weak solutions to the drift--diffusion system have been studied in \cite{mock_1974_initial, gajewski_1985_existence,gajewski_groger_1986}.
    The large time behavior of the isothermal drift--diffusion system \eqref{eq:NPV} has been studied in~\cite{gajewski_1996_reaction}. It has been proved that the solution to the transient system converges to the thermal equilibrium state as $t$  goes to infinity if the boundary conditions $N^D$, $P^D$, $V^D$ are in thermal equilibrium. 
    
   More precisely, the thermal equilibrium $(N^{eq}, P^{eq}, V^{eq})$ is a particular steady--state of \eqref{eq:NPV} for which electron and hole currents vanish, namely
    \[\nabla N^{eq}-N^{eq}\nabla V^{eq}\ =\ \nabla P^{eq}+P^{eq}\nabla V^{eq}\ =\ 0\,.\]

    If the Dirichlet boundary conditions satisfy $N^{D}$, $P^{D}>0$ and the compatibility conditions
    \begin{equation}\label{eq:hyp-compatibility-BC}
      \log(N^{D})-V^{D}=\alpha_{N} \text{ and } \log(P^{D})+V^{D}=\alpha_{P} \text{ on }\Gamma^{D},
    \end{equation}
    then the thermal equilibrium is defined as the solution to
    \begin{gather}
      -\lambda^{2}\Delta V^{eq}\ =\ \exp(\alpha_{P}-V^{eq})\,-\,\exp(\alpha_{N}+V^{eq})\,+\,C,\label{eqtherm-V}\\
      N^{eq}\ =\ \exp(\alpha_{N}+V^{eq}),\label{eqtherm-N}\\
      P^{eq}\ =\ \exp(\alpha_{P}-V^{eq}),\label{eqtherm-P}
    \end{gather}
    with Dirichlet-Neumann boundary conditions $V=V^D$ on $\Gamma_D$ and $\nabla V\cdot {\mathbf n}=0$ on $\Gamma_N$. The existence of a thermal equilibrium has been studied in \cite{markowich_1986_stationary, markowich_1990_semiconductor}.

    The proof of convergence to the thermal equilibrium is based on a relative entropy method, described for instance in the review paper \cite{Arnold2004}. Here the functional reads
    \begin{multline}\label{defEcontinu}
      \mathbb{E}(t)=\int_{\Omega}\left(\vphantom{\frac{1}{2}}H(N(t))-H(N^{eq})-\log(N^{eq})(N(t)-N^{eq})\right.\\
      + H(P(t))-H(P^{eq})-\log(P^{eq})(P(t)-P^{eq})\\
      \left.+\frac{\lambda^{2}}{2}|\nabla(V(t)-V^{eq})|^{2}\right)\dx\,,
    \end{multline}
    with $H(x)=\ds\int_{1}^{x}\log(s)\,ds$, and the entropy production functional is given by
    \[
      \mathbb{I}(t)=\int_{\Omega}\left(N\,|\nabla(\log(N)-V)|^{2}\,+\,P\,|\nabla(\log(P)+V)|^{2}\right)\dx\,.
    \]
    The entropy--entropy production inequality writes
    \begin{equation}\label{ineg_EI_continu_1}
      0 \leq \mathbb{E}(t)+\int_{0}^{t}\mathbb{I}(s)\,ds\leq\mathbb{E}(0).
    \end{equation}
    Using the Poincaré inequality and uniform positive upper and lower bounds on the densities one can show that there is a constant $C_{EI}>0$ such that for all $t\geq0$ one has 
    \begin{equation}\label{ineg_EI_continu_2}
      \mathbb{E}(t)\leq C_{EI}\, \mathbb{I}(t)\,.
    \end{equation}
    The combination of \eqref{ineg_EI_continu_1} and \eqref{ineg_EI_continu_2} yields exponential decay towards the thermal equilibrium. We refer to \cite{gajewski_groger_1986} for more details.

    \subsection{Numerical scheme}
    
    Let us define at each time step $n\in\NN$ the approximate solution $(f_{K}^{n})_{K\in\T}$ for $f=N,\,P,\,V$ and the approximate values at the boundary $(f_{\sigma})_{\sigma\in\Ee^{D}}$ (which do not depend on $n$ since the boundary data do not depend on time).  First of all, we discretize the initial and boundary conditions: 
    \begin{align}
      &\left(N_{K}^{0},P_{K}^{0}\right)\ =\ \frac{1}{m(K)}\int_{K}\left(N_{0},P_{0}\right)\,,\quad \forall K \in\T, \nonumber\\
      &\left(N_{\sigma}^D,P_{\sigma}^D,V_{\sigma}^D\right)\ =\ \frac{1}{m(\sigma)}\int_{\sigma}\left(N^{D},P^{D},V^{D}\right)\,,\quad\forall\sigma\in\Ee^{D}\,. \nonumber
    \end{align}
    Then, as for the previous models in Section~\ref{sec:fpe} and \ref{sec:pme}, we consider a backward Euler in time and finite volume in space discretization of the drift--diffusion system \eqref{eq:NPV}. The scheme writes
    \begin{align}
      &m(K) \frac{N_K^{n+1}-N_K^n}{\Delta t}\,+\,\ds\sum_{\sigma\in {\mathcal E}_K}{\mathcal F}_{K,\sigma}^{n+1}\ =\ 0\,,\ \forall K\in\T, \forall n\geq 0,\label{scheme-N}\\
      &m(K) \frac{P_K^{n+1}-P_K^n}{\Delta t}\,+\,\ds\sum_{\sigma\in {\mathcal E}_K}{\mathcal G}_{K,\sigma}^{n+1}\ =\ 0\,,\ \forall K\in\T, \forall n\geq 0,\label{scheme-P}\\
      &-\lambda^2\sum_{\sigma\in {\mathcal E}_K}\tau_{\sigma} D_{K,\sigma}V^{n}\ =\ m(K) (P_K^{n}-N_K^{n}+C_{K}),\ \forall K\in\T,\forall n\geq 0.\label{scheme-V}
    \end{align}
 For all $K\in\T$ and $\sigma\in\E[K]$, the numerical fluxes are defined by
    \begin{align}
      &\label{FLUX-SG-N}
      {\mathcal F}_{K,\sigma}^{n+1}=
      \tau_{\sigma}\left[ B\left(-D_{K,\sigma}V^{n+1}\right)N_K^{n+1}-B\left(D_{K,\sigma}V^{n+1}\right)N_{K,\sigma}^{n+1} \right],\\
      &\label{FLUX-SG-P}
      {\mathcal G}_{K,\sigma}^{n+1}=
      \tau_{\sigma}\left[ B\left(D_{K,\sigma}V^{n+1}\right)P_K^{n+1}-B\left(-D_{K,\sigma}V^{n+1}\right)P_{K,\sigma}^{n+1} \right].
    \end{align}
 The notations $f_{K,\sigma}^n$ and $D_{K,\sigma}f^n$ for $f=N,P,V $ are defined in \eqref{e:discretebc} and \eqref{e:diffop}. The function $B$ is any real function satisfying the properties \eqref{h:hyp_B}. As in Section~\ref{sec:fpe}, this notation encapsulates classical schemes such as the upwind, centered and Scharfetter-Gummel (SG) flux discretizations (see \eqref{e:Bfunctions}). 
    
    As in the previous sections we introduce the steady version of the numerical scheme \eqref{scheme-N}-\eqref{scheme-V}. It reads,
    \begin{align}
      &\ds\sum_{\sigma\in {\mathcal E}_K}{\mathcal F}_{K,\sigma}^{\infty}\ =\ 0\,,\label{scheme-steady-N}\\
      &\ds\sum_{\sigma\in {\mathcal E}_K}{\mathcal G}_{K,\sigma}^{\infty}\ =\ 0\,,\label{scheme-steady-P}\\
      &-\lambda^2\sum_{\sigma\in {\mathcal E}_K}\tau_{\sigma} D_{K,\sigma}V^{\infty}\ =\ m(K) (P_K^{\infty}-N_K^{\infty}+C_{K})\,,\label{scheme-steady-V}
    \end{align}
    with the corresponding fluxes
    \begin{align}
      &\label{FLUX-SG-N-steady}
      {\mathcal F}_{K,\sigma}^{\infty}=
      \tau_{\sigma}\left[ B\left(-D_{K,\sigma}V^{\infty}\right)N_K^{\infty}-B\left(D_{K,\sigma}V^{\infty}\right)N_{K,\sigma}^{\infty} \right],\\
      &\label{FLUX-SG-P-steady}
      {\mathcal G}_{K,\sigma}^{\infty}=
      \tau_{\sigma}\left[ B\left(D_{K,\sigma}V^{\infty}\right)P_K^{\infty}-B\left(-D_{K,\sigma}V^{\infty}\right)P_{K,\sigma}^{\infty} \right].
    \end{align}

Let us just recall that the fundamental property of the Scharfetter-Gummel scheme, with $B(x) = x/(\exp(x)-1)$ is that for thermal boundary conditions \eqref{eq:hyp-compatibility-BC} the discrete steady state solving \eqref{scheme-steady-N}--\eqref{scheme-steady-V} is of the form $(N^\infty, P^\infty, V^\infty) = (N^\text{eq}, P^\text{eq}, V^\text{eq})$ with 
    \be\label{discr_thermal_PN}
    N_K^\text{eq}\ =\ \exp(\alpha_{N}+V^\text{eq}_K)\text{ and }P_K^\text{eq}\ =\ \exp(\alpha_{P}-V^\text{eq}_K)\,,\ \forall K\in\T.
    \ee
    In other words, the equivalents of \eqref{eqtherm-N} and \eqref{eqtherm-P} hold at the discrete level. Therefore, in the case of the SG scheme the resolution of the discrete steady system \eqref{scheme-steady-N}-\eqref{scheme-steady-V} amounts to solving the nonlinear system of equations
    \be\label{discr_thermal_V}
    -\lambda^2\sum_{\sigma\in {\mathcal E}_K}\tau_{\sigma} D_{K,\sigma}V^\text{eq}\ =\ m(K) (\exp(\alpha_{P}-V^\text{eq}_K)-\exp(\alpha_{N}+V^\text{eq}_K)+C_{K}),\ \forall K\in\T,
    \ee
 which is the discrete counterpart of \eqref{eqtherm-V}. For general $B$-functions, the discrete steady state $(N^\infty, P^\infty, V^\infty)$ differs from the discrete thermal equilibrium $(N^\text{eq}, P^\text{eq}, V^\text{eq})$.

    Another important fact is that the discrete equivalent of the entropy--entropy production principle \eqref{ineg_EI_continu_1}  as well as the exponential decay of the discrete relative entropy were only proved and illustrated numerically in the context of SG schemes, see \cite{chatard_fvca6,bessemoulin_chainais_2017}. For $n\in\NN$, we may define two discrete relative entropy functional $ E_\text{eq}^{n}$ and $  E_\infty^{n}$ depending on the choice of stationary state. They read, with $\text{rel}=\text{eq}$ or $\infty$, 
    \begin{multline}\label{defEdiscret}
      E_\text{rel}^{n}\ =\ \sum_{K\in\T}m(K)\left[\vphantom{\frac12}H(N_K^n)-H(N^{\text{rel}}_K)-\log(N^{\text{rel}}_K)(N_K^n-N^{\text{rel}}_K)\right.\\
      + \left.\vphantom{\frac12}H(P_K^n)-H(P^{\text{rel}}_K)-\log(P^{\text{rel}}_K)(P_K^n-P^{\text{rel}}_K)\right]\\
      +\frac{\lambda^{2}}{2}\sum_{\sigma\in\E}\tau_\sigma\,|D_{K,\sigma}(V^n-V^{\text{rel}})|^{2}\,,
    \end{multline}
    In the case of the Scharfetter-Gummel scheme and when the boundary data satisfy the thermal equilibrium,  one has $E_\infty^{n} = E_{eq}^{n}$. Under these assumptions, the exponential decay of the relative discrete entropy has been proved in \cite{bessemoulin_chainais_2017}, namely 
      \[
          E_\text{eq}^{n}\ \leq\ E_\text{eq}^{0}\,e^{-\kappa\,t^n}\,,
    \]
   with a rate $\kappa$ independent of the size of the discretization. As a consequence the discrete densities and potential converge to the discrete thermal equilibrium at exponential rate.

    \medskip
    
    In the last part of this paper, we provide numerical evidence that exponential decay of $E_\infty^n$ also occurs for the other monotone numerical schemes such as the upwind or centered schemes. Let us emphasize once again that for  other $B$-schemes we consider  the entropy $E_\infty^n$ relative to the steady state of the scheme, which solves \eqref{scheme-steady-N}--\eqref{scheme-steady-V}. Except for the Scharfetter-Gummel scheme, it does not coincide with $E^n_{eq}$.
    
    As in Section~\ref{sec:fpe}, it emphasizes that the accurate approximation of the discrete steady state and exponential return to equilibrium are not the same notions. Of course, the upwind and centered schemes do not enjoy the property of preserving the thermal equilibrium relation  \eqref{eqtherm-N} and \eqref{eqtherm-P} as the SG scheme does. However, the monotonous structure of these schemes seems to be sufficient to get exponential decay rate of the relative entropy without saturation before machine precision.

    \subsection{Numerical results}
    
   The schemes are nonlinear and we proceed with a Newton-Raphson method to solve the nonlinear system of equations at each time step. The stationary system is also solved thanks to a Newton-Raphson method. 
   
   In every test cases we choose a fixed time step $\Delta t = 10^{-2}$.  The domain and the mesh structure are the same as in Section~\ref{sec:numresFP}. The mesh is taken such that $\Delta x = 1/32$ which corresponds to $3584$ triangles.
    
    \subsubsection{Exponential decay for various $B$-schemes}
    
    Our first test case is taken from \cite{bessemoulin_chainais_2017}. 
  
         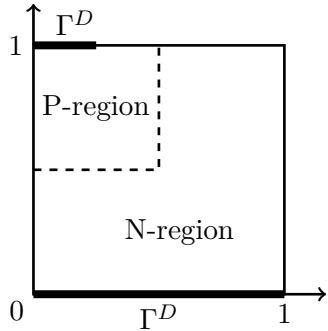
\begin{figure}[ht!]  
    \begin{center}
\begin{tikzpicture}[scale=1.1]
\draw[line width=1pt,->] (0.,3)--(0.,3.5);
\draw[line width=1pt,->] (3,0)--(3.5,0.);
\draw[line width=1pt] (0,0) rectangle (3,3);
\draw[line width=1pt,dashed] (0,1.5)--(1.5,1.5)--(1.5,3); 
\draw[line width=3pt] (0,3)--(0.75,3);
\draw[line width=3pt] (0,0)--(3,0);
\node at (0.75,2.25){P-region};
\node at (1.75,.75){N-region};
\node[above] at (0.5,3){$\Gamma^D$};
\node[below] at (1.5,0){$\Gamma^D$};
\node at (-0.2,-0.2){0};
\node[below] at (3,0){1};
\node[left] at (0,3){1};
\end{tikzpicture}
\end{center}
     \caption{The PN junction diode}\label{f:PNJunc}
        \end{figure}

\begin{figure}[ht!]
     \begin{tabular}{ccc}
      \includegraphics[width = .33\linewidth]{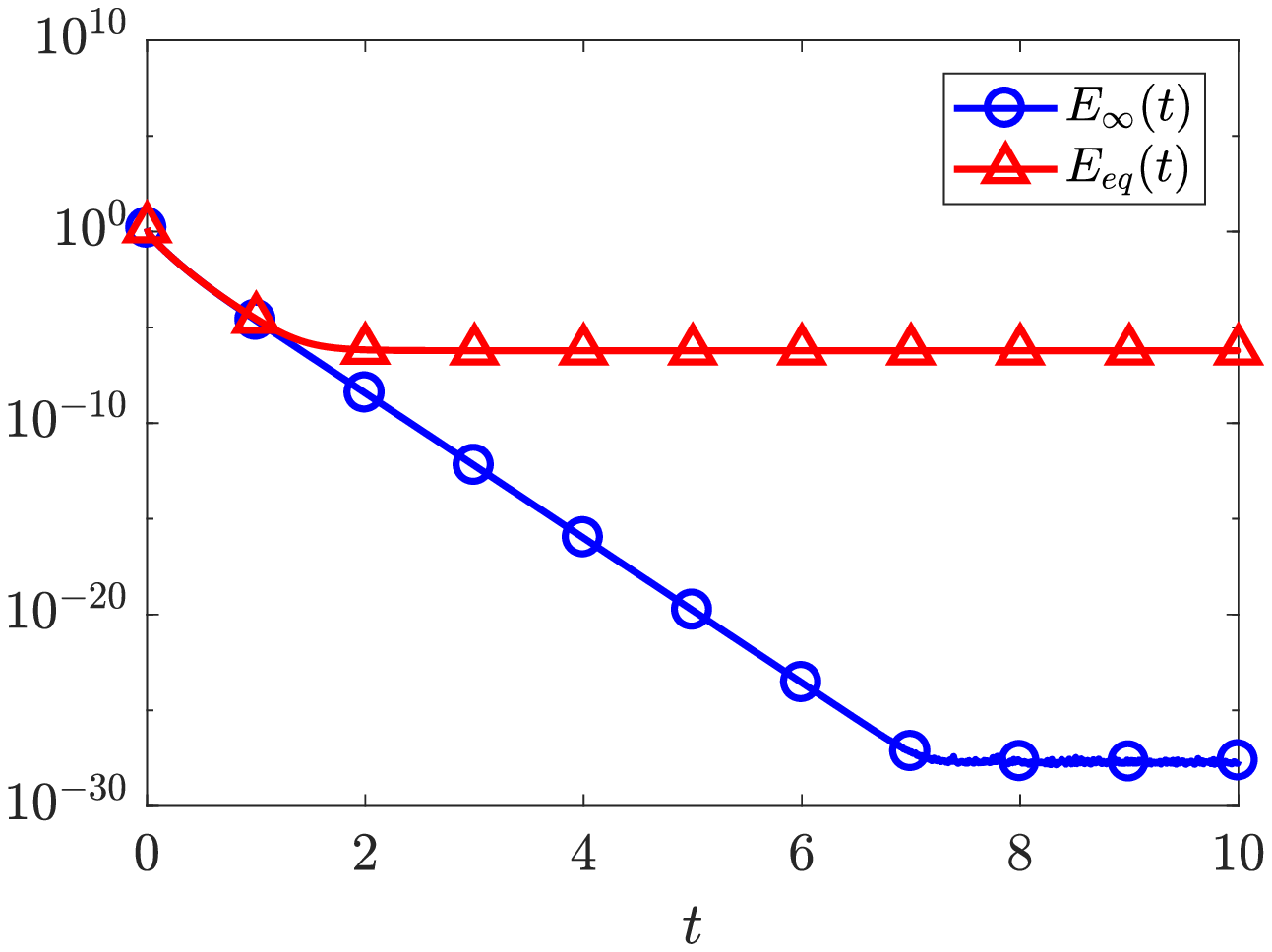}& 
      \includegraphics[width = .33\linewidth]{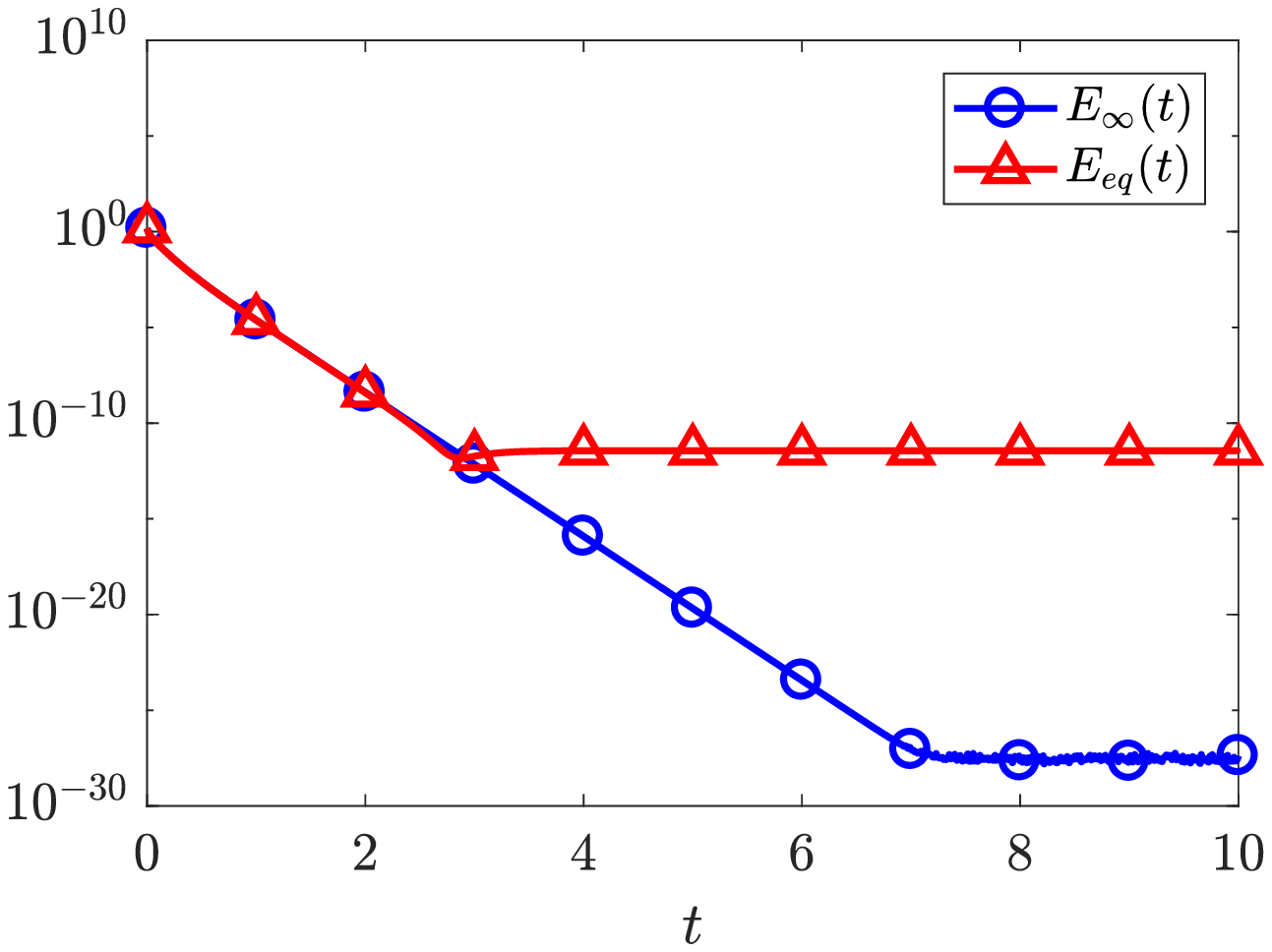}&
      \includegraphics[width = .33\linewidth]{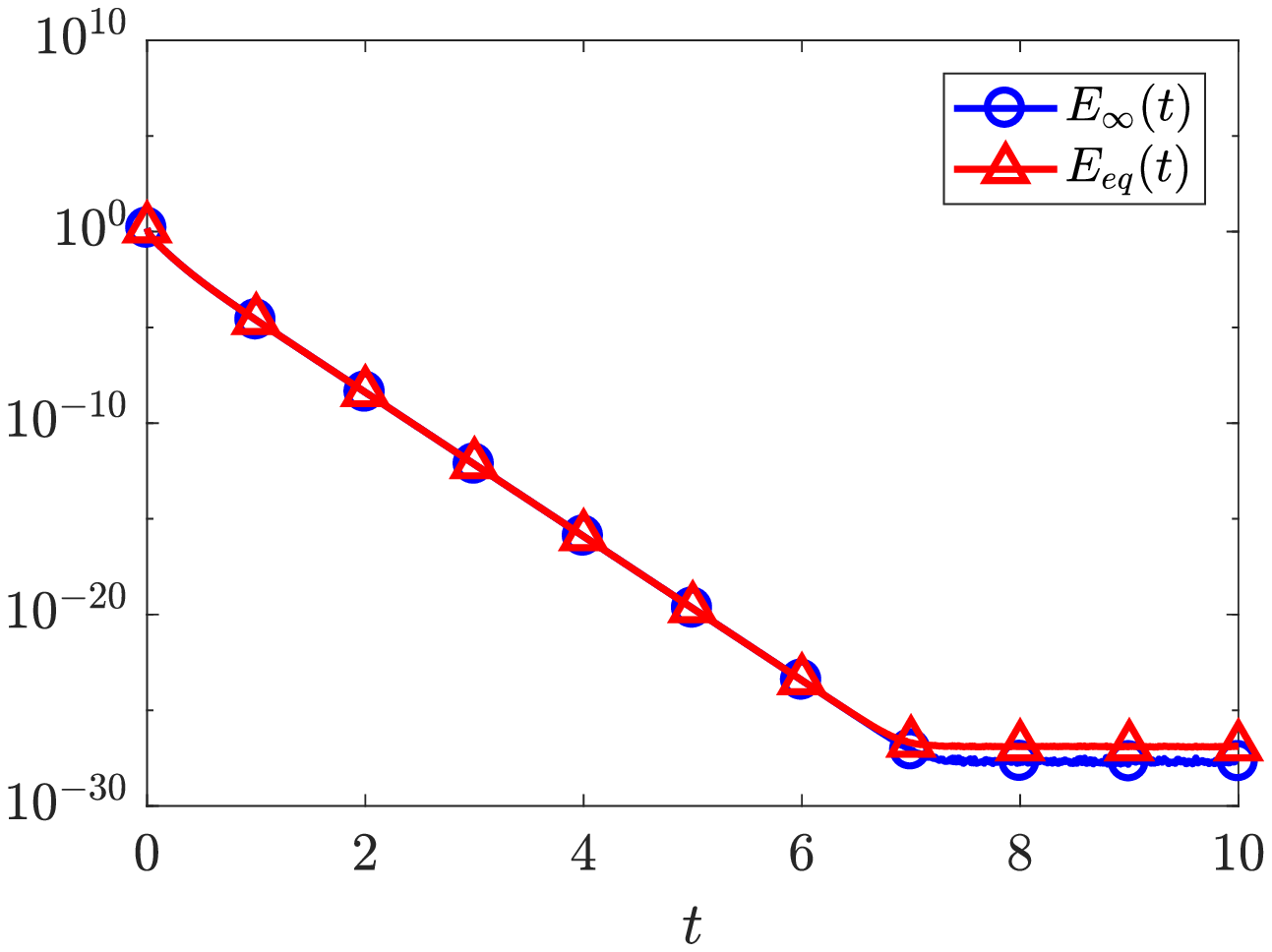}\\
      (a) Upwind&(b) Centered&(c) Scharfetter-Gummel
     \end{tabular}
     \caption{Evolution of the discrete entropy functional versus time for various $B$-schemes. The entropies are relative either to the steady state of the scheme ($E_\infty(t)$) or to the discrete thermal equilibrium ($E_{eq}(t)$).}\label{f:dd_entropy}
    \end{figure}

    The model is a PN-junction in 2D as shown on Figure~\ref{f:PNJunc} with Dirchlet boundary conditions on $\Gamma_D = \{x_2 = 0\}\cup\{x_2=1, 0\leq x_1\leq 0.25\}$. The boundary conditions on the electron density are taken with $N^D(x_1,x_2)= e$ and $P^D(x_1,x_2)= e^{-1}$ if $x_2 = 0$ and $N^D(x_1,x_2)= P^D(x_1,x_2) = 1$ if $x_2 = 1$  and $0\leq x_1\leq 0.25$. The boundary condition on the potential is taken as $V^D = (\log(N^D) - \log(P^D))/2$. Observe that the compatibility condition \eqref{eq:hyp-compatibility-BC} is satisfied with $\alpha_N = \alpha_P = 0$. The initial data on the densities are $N_0(x_1,x_2) = e + (1-e)(1-\sqrt{x_2})$ and $P_0(x_1,x_2) = e^{-1} + (1-e^{-1})(1-\sqrt{x_2})$. On Figure~\ref{f:dd_entropy}, we show the evolution of the two kinds of relative entropies for each scheme. It seems that as announced the entropy relative to the discrete steady states always decays exponentially until machine precision. However the entropy relative to the thermal equilibrium saturates at some point, except in the case of the SG scheme for which both discrete steady states coincide.

    \subsubsection{Non-thermal steady states}
     \begin{figure}[ht!]
    \includegraphics[width = .5\linewidth]{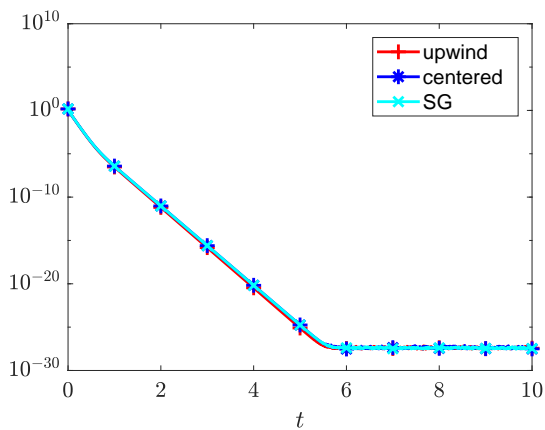}
     \caption{Evolution of the discrete relative entropy ${\mathbb E}_\infty$ versus time for various $B$-schemes.  The steady state is not a thermal equilibrium here. The entropy is relative to steady state of each scheme ($E_\infty$).}\label{f:dd_nonthermal}
    \end{figure}

    In this last test case we explore the behavior of the schemes in an out-of-equilibrium situation. More precisely, we consider the same test case as in the last paragraph except that we add a bias to the Dirichlet condition on the potential, $V^D(x_1,x_2) = (\log(N^D(x_1,x_2)) - \log(P^D(x_1,x_2)))/2 + V_{\text{bias}}(x_1,x_2)$, with $V_{\text{bias}}(x_1,x_2) = 2.5$ if $x_2 = 0$ and $V_{\text{bias}}(x_1,x_2) = -2.5$ if $x_2 = 1$. In this context, the compatibility condition \eqref{eq:hyp-compatibility-BC} is not satisfied anymore and the steady state is not a thermal equilibrium. At the continuous level one can show existence of such steady states (see \cite{markowich_1993_vacuum}), however uniqueness may be lost for large bias potential (see for instance \cite{mock_1982_example}). To our knowledge, there are no theoretical results concerning the large time behavior in non-thermal situations. On Figure~\ref{f:dd_nonthermal}, our numerical simulations show that, for this particular test case, exponential decay still seems to occur.
    
    \section{Conclusion}\label{sec:conclusion}
    
    In this paper, we have considered several convection-diffusion models set on a bounded domain with mixed non-homogeneous Dirichlet and no-flux boundary conditions. These various models included linear Fokker-Planck equations, non-linear porous media equations, and nonlinear systems of drift-diffusion equations. They all satisfied a relative entropy principle at the continuous level, which allows to prove and quantify the exponential return to a steady-state in the large.
    
    We have considered a class of numerical schemes for these equations which was of finite volume type and implicit in time. The numerical fluxes were discretized with two-point monotone fluxes written in the general framework of $B$-schemes, which allows to deal with upwind, centered and Schartfetter-Gummel discretizations.
    
    For the linear and non-linear equations we have shown theoretically that the relative entropy principles also hold at the discrete level and that the discrete solution to all these schemes exponentially to the discrete steady state of the scheme. For the drift-diffusion system, we illustrated the same phenomena with numerical simulations. 
    
    This contribution primarily emphasized that the property of exponential decay to the steady state may be obtained for a large class of schemes. The particular choice of scheme (such as Scharfetter-Gummel for drift-diffusion) mainly provides a better approximation of the steady state. Besides, we provided new theoretical results on the large-time behavior of these classical schemes when used for boundary-driven convection-diffusions.  For Fokker-Planck and porous media equations, we were able to derive explicit discrete estimates for decay rates depending on parameters of the equations and bounds on the boundary data.

    \bibliographystyle{myplain}
    \bibliography{bibli}
  \end{document}